\documentclass[reqno,a4paper,11pt]{amsart}

\usepackage[utf8]{inputenc}
\usepackage[T1]{fontenc}
\usepackage[english]{babel}
\usepackage[babel]{microtype}
\usepackage[a4paper, margin=1in]{geometry}
\usepackage[dvipsnames]{xcolor}

\usepackage{csquotes}

\usepackage{amsmath,amssymb,amsthm}

\usepackage{mathrsfs}

\usepackage{mathtools}
\usepackage{commath}

\usepackage{thmtools}
\usepackage{thm-restate}

\usepackage{bbm}
\usepackage{dsfont}
\usepackage{lmodern}
\usepackage{fixcmex}

\usepackage{cases}
\usepackage{enumitem}
\setlist[enumerate,1]{label=(\roman*)}

\usepackage{centernot}
\usepackage{booktabs}
\usepackage{multirow}
\usepackage{comment}
\usepackage{float}
\usepackage{scalerel}

\usepackage[pdfborderstyle={/S/U/W 0},unicode]{hyperref}
\hypersetup{
    colorlinks=true,
    linkcolor=blue!85!black,
    citecolor=blue!85!black,
    urlcolor=blue!85!black,
}

\usepackage{cleveref}

\usepackage{etoolbox}

\usepackage{subfiles}
\usepackage{pdfpages}

\numberwithin{equation}{section}

\ifdefined\checkunusedreferences
\usepackage{refcheck}

\makeatletter
\newcommand{\alsocheck}[1]{%
  \expandafter\let\csname @@\string#1\endcsname#1%
  \expandafter\DeclareRobustCommand\csname relax\string#1\endcsname[1]{%
    \csname @@\string#1\endcsname{##1}\@for\@temp:=##1\do{\wrtusdrf{\@temp}\wrtusdrf{{\@temp}}}}%
  \expandafter\let\expandafter#1\csname relax\string#1\endcsname
}
\newcommand{\alsocheckrange}[1]{%
  \expandafter\let\csname @@\string#1\endcsname#1%
  \expandafter\DeclareRobustCommand\csname relax\string#1\endcsname[2]{%
    \csname @@\string#1\endcsname{##1}{##2}\wrtusdrf{##1}\wrtusdrf{{##1}}\wrtusdrf{##2}\wrtusdrf{{##2}}}%
  \expandafter\let\expandafter#1\csname relax\string#1\endcsname
}
\makeatother

\alsocheck{\cref}
\alsocheck{\Cref}
\alsocheckrange{\crefrange}
\alsocheckrange{\Crefrange}

\else
\fi

\ifdefined\defcolor
\declaretheoremstyle[
  headfont=\normalfont\bfseries,
  bodyfont=\normalfont,
  shaded={bgcolor=\defcolor}
]{noital}
\else
\declaretheoremstyle[
  headfont=\normalfont\bfseries,
  bodyfont=\normalfont,
]{noital}
\fi

\declaretheorem[style=plain,numberwithin=section,name=Theorem]{theorem}

\declaretheorem[style=plain,sibling=theorem,name=Proposition]{proposition}

\declaretheorem[style=plain,sibling=theorem,name=Lemma]{lemma}

\declaretheorem[style=plain,sibling=theorem,name=Question]{question}

\declaretheorem[style=plain,sibling=theorem,name=Claim]{claim}

\declaretheorem[style=plain,numbered=no,name=Theorem]{theoremn-n}
\declaretheorem[style=plain,numbered=no,name=Theorems]{theorems-n}
\declaretheorem[style=plain,numbered=no,name=Proposition]{proposition-n}
\declaretheorem[style=plain,numbered=no,name=Propositions]{propositions-n}
\declaretheorem[style=plain,numbered=no,name=Lemma]{lemma-n}
\declaretheorem[style=plain,numbered=no,name=Lemmas]{lemmas-n}
\declaretheorem[style=plain,numbered=no,name=Corollary]{corollary-n}
\declaretheorem[style=plain,numbered=no,name=Corollaries]{corollaries-n}
\declaretheorem[style=plain,numbered=no,name=Conjecture]{conjecture-n}
\declaretheorem[style=plain,numbered=no,name=Conjectures]{conjectures-n}
\declaretheorem[style=plain,numbered=no,name=Question]{question-n}
\declaretheorem[style=plain,numbered=no,name=Questions]{questions-n}
\declaretheorem[style=plain,numbered=no,name=Claim]{claim-n}
\declaretheorem[style=plain,numbered=no,name=Claims]{claims-n}
\declaretheorem[style=plain,numbered=no,name=Fact]{fact-n}
\declaretheorem[style=plain,numbered=no,name=Facts]{facts-n}
\declaretheorem[style=plain,numbered=no,name=Problem]{problem-n}
\declaretheorem[style=plain,numbered=no,name=Problems]{problems-n}
\declaretheorem[style=plain,numbered=no,name=Open Problem]{openproblem-n}
\declaretheorem[style=plain,numbered=no,name=Open Problems]{openproblems-n}
\declaretheorem[style=plain,numbered=no,name=Challenge]{challenge-n}
\declaretheorem[style=plain,numbered=no,name=Challenges]{challenges-n}
\declaretheorem[style=plain,numbered=no,name=Exercise]{exercise-n}
\declaretheorem[style=plain,numbered=no,name=Exercises]{exercises-n}
\declaretheorem[style=plain,numbered=no,name=Property]{property-n}
\declaretheorem[style=plain,numbered=no,name=Properties]{properties-n}

\declaretheorem[style=noital,numbered=no,name=Remark]{remark-n}
\declaretheorem[style=noital,numbered=no,name=Remarks]{remarks-n}
\declaretheorem[style=noital,numbered=no,name=Definition]{definition-n}
\declaretheorem[style=noital,numbered=no,name=Definitions]{definitions-n}
\declaretheorem[style=noital,numbered=no,name=Construction]{construction-n}
\declaretheorem[style=noital,numbered=no,name=Constructions]{constructions-n}
\declaretheorem[style=noital,numbered=no,name=Observation]{observation-n}
\declaretheorem[style=noital,numbered=no,name=Observations]{observations-n}
\declaretheorem[style=noital,numbered=no,name=Example]{example-n}
\declaretheorem[style=noital,numbered=no,name=Examples]{examples-n}

\renewcommand{\le}{\leqslant}

\renewcommand{\ge}{\geqslant}

\let\oldexists\exists
\let\exists\relax
\DeclareMathOperator{\exists}{\:\!\oldexists}
\let\oldforall\forall
\let\forall\relax
\DeclareMathOperator{\forall}{\:\!\oldforall}

\undef{\set}
\DeclarePairedDelimiter{\set}{\lbrace}{\rbrace}

\undef{\emptyset}
\newcommand{\emptyset}{\varnothing}

\DeclareMathOperator{\ind}{\mathbf{1}}

\renewcommand{\d}{\mathop{}\!\mathrm{d}}

\undef{\mod}
\newcommand{\mod}[1]{\ (\mathrm{mod}\ #1)}

\undef{\abs}
\DeclarePairedDelimiterX{\abs}[1]
  {\lvert}{\rvert}{\ifblank{#1}{\,\cdot\,}{#1}}
\undef{\norm}
\DeclarePairedDelimiterX{\norm}[1]
  {\lVert}{\rVert}{\ifblank{#1}{\,\cdot\,}{#1}}
\DeclarePairedDelimiterX{\inner}[2]
  {\langle}{\rangle}{\ifblank{#1}{\,\cdot\,}{#1},\ifblank{#2}{\,\cdot\,}{#2}}

\DeclareMathOperator{\diam}{diam}

\DeclareMathDelimiter{\given}
  {\mathbin}{symbols}{"6A}{largesymbols}{"0C}
\DeclareMathOperator{\Prob}{\mathbb{P}}
\DeclarePairedDelimiterXPP{\prob}[1]
  {\Prob}{\lparen}{\rparen}{}
  {\renewcommand{\given}{\nonscript\;\delimsize\vert\nonscript\;\mathopen{}}#1}
\DeclareMathOperator{\Expec}{\mathbb{E}}
\DeclarePairedDelimiterXPP{\expec}[1]
  {\Expec}{\lparen}{\rparen}{}
  {\renewcommand{\given}{\nonscript\;\delimsize\vert\nonscript\;\mathopen{}}#1}
\DeclareMathOperator{\Var}{Var}
\DeclarePairedDelimiterXPP{\var}[1]
  {\Var}{\lparen}{\rparen}{}
  {\renewcommand{\given}{\nonscript\;\delimsize\vert\nonscript\;\mathopen{}}#1}
\DeclareMathOperator{\Cov}{Cov}
\DeclarePairedDelimiterXPP{\cov}[2]
  {\Cov}{\lparen}{\rparen}{}{#1,#2}

\newcommand{\eps}{\varepsilon}

\newcommand{\EE}{\mathbb{E}}

\newcommand{\NN}{\mathbb{N}}

\newcommand{\PP}{\mathbb{P}}

\newcommand{\RR}{\mathbb{R}}

\newcommand{\ZZ}{\mathbb{Z}}

\newcommand{\cA}{\mathcal{A}}

\newcommand{\cC}{\mathcal{C}}

\newcommand{\cE}{\mathcal{E}}

\newcommand{\cM}{\mathcal{M}}
\newcommand{\cN}{\mathcal{N}}

\newcommand{\bfe}{\mathbf{e}}

\newcommand{\bfu}{\mathbf{u}}
\newcommand{\bfv}{\mathbf{v}}
\newcommand{\bfw}{\mathbf{w}}

\newcommand{\cgH}{c(\mathbb{H})}
\newcommand{\gH}{\mathbb{H}}



\usepackage{tikz}
\usetikzlibrary{arrows, arrows.meta, backgrounds}

\DeclareMathOperator{\Diam}{Diam}

\renewcommand{\complement}{\mathsf{c}}
\newcommand{\e}{\textrm{e}}

\newcommand{\symdiff}{\mathbin{\triangle}}

\begin{document}

\title[Source localisation in simple random walks]{Source localisation in simple random walks}

\author[R. Goenka]{Ritesh Goenka}
\author[P. Keevash]{Peter Keevash}
\author[T. Przybyłowski]{Tomasz Przybyłowski}

\address{Mathematical Institute, University of Oxford, Oxford, United Kingdom}
\email{\{goenka,keevash,przybylowski\}@maths.ox.ac.uk}

\begin{abstract}
    We consider the problem of locating the source (starting vertex) of a simple random walk, given a snapshot of the set of edges (or vertices) visited in the first $n$ steps. Considering lattices $\ZZ^d$, in dimensions $d \ge 5$, we show that the source can be identified (a) with probability bounded away from $0$ using one guess, and (b) with probability arbitrarily close to $1$ using a constant number of guesses. On the other hand, for dimensions $d \le 2$, we show that one cannot locate the source with positive constant probability. Our arguments apply more generally to strongly transient and recurrent simple random walks on vertex-transitive graphs.
\end{abstract}

\subjclass[2020]{Primary: 60G50; Secondary: 60C05}
\keywords{Random walk, source detection, trace, graph, cut-edges.}

\maketitle

\section{Introduction}
\label{sec:intro}

Growth processes are ubiquitous in nature, with examples ranging from the growth of cancer tumours in a lab to polymer formation~\cite{Tur}. These phenomena have inspired many mathematical models of random growth processes, including Diffusion Limited Aggregation (DLA) \cite{WS}, First Passage Percolation (FPP) \cite{HWe}, Dielectric Breakdown Model (DBM) \cite{NPW}, Eden growth process \cite{Eden}, Internal DLA (IDLA) \cite{LBG}, and the Simple Random Walk (SRW) viewed at first passage times. Depending on the underlying graph $\gH$ and the stochastic growth rule, the growth cluster may exhibit very different kinds of typical structure. In particular, on lattices $\gH = \ZZ^d$, the growth clusters in DLA, IDLA, and the high-dimensional simple random walk viewed at first passage times have fractal, Euclidean ball, and quasi-linear structures, respectively.

Our focus in this paper is the question of locating the source (starting vertex) of a growth process given a snapshot of the cluster at some large (or infinite) time. For a random process, this is naturally a statistical question, in which one aims for a prediction that is successful with constant or high probability for cluster size $n \to \infty$. There is an extensive empirical literature (see the survey \cite{SaCh}) in the context of rumor source detection on various real-world networks and also theoretical analysis of some models, such as the Eden growth process on the infinite $d$-regular tree~\cite{ShZ,ShZ2} and an extensive literature on finding the root of randomly growing trees; see~\cite{BDL17,BH23,BB22,CCL+24}. We consider this question for simple random walks.

\subsection{Localisation for the SRW}
Throughout the paper, we consider the simple random walk $(X_n)_{n \ge 0}$ on an infinite vertex-transitive graph $\gH$ of finite degree. This walk is called  \emph{recurrent} if it returns to the source almost surely, or \emph{transient} otherwise. It is called \emph{strongly transient} if it is transient and the expected return time, given that it returns, is finite; for example, lattices $\ZZ^d$ with $d \ge 5$ are strongly transient (see \Cref{app:strongtrans} for details). Let $G_n = G^{\gH}_n$ be the $n$-step \emph{trace} of the walk, namely, the graph that consists of vertices and edges visited by the walk in the first $n$ steps, embedded inside the host graph $\gH$. We write $\PP_v$ for the law of the walk started at a vertex $X_0 = v \in \gH$, which we call the \emph{source}. An \emph{estimator} for the source is a function $\Phi(G)$ that takes an embedded graph $G$ as input and outputs a vertex (or a set of vertices) in $\gH$, which is possibly random conditional on the trace. We evaluate the performance of $\Phi$ by $ \inf_{v \in V(\gH)} \PP_v(\Phi(G_n) = v)$. If $\Phi$ is translation invariant, we may simply write $\PP_v(\Phi(G_n) = v)$.

Our main result shows that if $\gH$ is strongly transient, then the source can be located with constant probability, using the following natural estimator $\Psi$: given an input graph $G$, select a uniformly random diametrical (maximum distance) path and output one of its endpoints uniformly at random. We note that this estimator only depends on $G$, not its embedding in $\gH$.

\begin{theorem}
\label{thm:SRW}
    For any infinite vertex-transitive strongly transient graph $\gH$ of finite degree,
    \begin{equation*}
        \liminf_{n \to \infty} \PP_v\big(\Psi(G_n^\gH) = v\big) > 0.
    \end{equation*}
    In particular, source localisation is possible if $\gH = \ZZ^d$ for $d \ge 5$.
\end{theorem}

For $\gH = \ZZ^d$, the estimator $\Psi$ yields a localisation probability of $1/2 - O(1/d)$ as $d \rightarrow \infty$; we show that this is sharp.

\begin{theorem}
\label{thm:SRWoptimal}
    We have
    \begin{equation*}
        \liminf_{n\to\infty} \inf_{v \in \ZZ^d} \PP_v\big(\Psi\big(G_n^{\ZZ^d}\big) = v\big) = 1/2 - O(1/d) \ \text{ as } \ d \to \infty.
    \end{equation*}
    Moreover, for every estimator $\Phi$, we have 
    \begin{equation*}
        \limsup_{n\to\infty} \inf_{v \in \ZZ^d} \PP_v(\Phi(G_n^{\ZZ^d}) = v) = 1/2 - \Omega(1/d) \ \text{ as } \ d \to \infty.
    \end{equation*}
\end{theorem}

Our next result concerns the impossibility of localising the source. We say that the SRW on $\gH$ is \emph{amnesiac} if $\limsup_{n \to \infty} \inf_{v \in V(\gH)} \PP_v(\Phi(G^\gH_n) = v) = 0$ for every estimator~$\Phi$.

\begin{theorem}
\label{thm:SRWrecurrent}
    Suppose $\gH$ is a recurrent, infinite, connected, vertex-transitive graph of finite degree. Then the simple random walk on $\gH$ is amnesiac.
\end{theorem}

In particular, the above theorem rules out localisation in $\ZZ^1$ and $\ZZ^2$, but leaves $\ZZ^3$ and $\ZZ^4$ unresolved (see \Cref{sec:conclude}).

\subsection{Overview of ideas and methods}
\label{sec:overview}

The starting point for constructing an estimator for locating the source in \Cref{thm:SRW} is the observation that by reversibility the source vertex $X_0$ and the final vertex $X_n$ are indistinguishable, so it is equivalent (and more convenient) to locate $X_n$ rather than $X_0$. The success of the estimator $\Psi$ based on diametrical paths may be intuitively understood from the phenomenon that in a strongly transient graph, the number of cut-edges of the trace $G_n$ grows linearly in $n$. This suggests that there should be a constant probability of the event $\cC_{n-1}$ that step $n$ increases the diameter of the trace, in which case $X_n$ is an endpoint of all diametrical paths.

However, there are significant challenges in making this heuristic rigorous. The most significant is that the only straightforward conclusion (via the ergodic theorem) is a lower bound for the Ces\`aro averages $\frac1n \sum_{i=0}^{n-1} \PP(\cC_i)$. The key technical ingredients of the proof concern the distribution of the set $L$ of cut edges of the trace, in particular, showing that with high probability, it has a large intersection with every diametrical path. This allows us to control the variation of $\PP(\cC_i)$ on linear scales, thus boosting Ces\`aro control to a lower bound on $\liminf_n \PP(\cC_n)$.

Our upper bounds on the localisation probability are based on coupling arguments. The simplest such approach is a global coupling based on reversibility, which implies that no estimator $\Phi$ for the SRW in $\ZZ^d$ can succeed with probability greater than $1/2$. We do not assume $\Phi$ to be translation invariant, which makes the argument slightly cumbersome. In our proof, we consider a large box $B \subseteq \ZZ^d$ and the random variable $W = \sum_{v \in B} [\mathbf{1}(\Phi(G_n + v) = v) + \mathbf{1}(\Phi(G'_n + v + X_n) = v + X_n)]$, where $G_n$ is the trace of $(X_0,\dots,X_n)$ and $G'_n$ is the trace of the reversed walk. Then $G_n + v$ and $G'_n + v + X_n$ are identical, so whenever $X_n \ne X_0$, we have $|W| \le |B|$. We then obtain the required bound from $\mathbb{E}W \ge (2+o(1)) |B| \inf_{v \in \ZZ^d} \PP(\Phi(G_n+v)=v)$ as $|B| \to \infty$; here the $o(1)$ term is justified by amenability of $\ZZ^d$. A more sophisticated coupling based on re-routing the ends of a walk leads to the improved bound required for \Cref{thm:SRWoptimal}.

For recurrent graphs, we obtain the bound in \Cref{thm:SRWrecurrent} by a similar calculation as above in combination with the idea that several walks starting at distinct vertices can, with high probability, be coupled so that they have identical traces.

\subsection{Further results}

Here we consider several natural variants of the source localisation problem for the simple random walk. For simplicity, we only consider $\gH = \ZZ^d$ for $d \ge 5$, but our arguments work for any strongly transient~$\gH$. For brevity, we shall write $G_n^d$ for $G_n^{\ZZ^d}$.

First we consider the \emph{high accuracy} problem, where given a target accuracy $\varepsilon > 0$, we construct an estimator $\Lambda_K(G_n^d)$ which is a set of $K = K(\varepsilon)$ vertices that contains the source with probability at least  $1 - \varepsilon$. Given an input graph $G$ and $k \ge 2$ even, we define $\Lambda_k(G)$ by selecting a uniformly random pair $(u,v)$ of diametric vertices in $G$ and outputting the $k/2$ closest vertices to $u$ and the $k/2$ closest vertices to $v$.

\begin{theorem}
\label{thm:SRWconst}
    Let $d \ge 5$. For each $\varepsilon \in (0, 1)$, there is a constant $K \coloneqq K(\varepsilon) \in \NN$ such that
    \begin{equation*}
        \liminf_{n \to \infty} \PP_v \big(v \in \Lambda_{K}(G_n^d)\big) \ge 1 - \varepsilon.
    \end{equation*}
\end{theorem}

Next we consider the problem of source localisation given the \emph{total} (infinite) trace $G_\infty^d$.

\begin{theorem}
\label{thm:infinite}
    Let $d \ge 5$. Then there is an estimator $\Gamma$ such that
    \begin{equation*}
        I_d \coloneqq \inf_{v \in \ZZ^d}  \PP_v(\Gamma(G_\infty^d) = v) \text{ satisfies } I_d  > 0 \text{ and } I_d = 1-O(1/d) \text{ as } d \to \infty.
    \end{equation*}
    Furthermore, $\inf_{v \in \ZZ^d} \PP_v(\Phi(G_\infty^d) = v) = 1 - \Omega(1/d)$ as $d \rightarrow \infty$, for any estimator $\Phi$.
\end{theorem}

We also show that the estimator $\Psi$ in \Cref{thm:SRW} is effective when given a trace with exactly $n$ vertices: formally $R_n^d \coloneqq G^d_{\tau_n}$ where $\tau_n = \inf \{t \ge 0 \colon |V(G_t^d)| = n\}$.

\begin{theorem}
\label{thm:range}
    Let $d \ge 5$. Then $\liminf_{n\to\infty} \PP_v(\Psi(R_n^d)=v) > 0$.
\end{theorem}

Our final variant concerns the harder source localisation problem where one cannot see the edges of the trace, only its vertex set $V_n$. It turns out that the same estimator $\Psi$ in \Cref{thm:SRW} applied to the induced subgraph $\widetilde{G}_n^d = \ZZ^d[V_n]$ is still effective!

\begin{theorem}
\label{thm:SRW2}
    Let $d \ge 5$. Then $\liminf_{n \to \infty} \PP_v\big(\Psi(\widetilde{G}_n^d) = v\big) > 0$.
\end{theorem}

In fact, \Cref{thm:infinite,thm:SRWconst,thm:range} also admit stronger versions with the vertex trace, which can be proved by the same proof technique as in \Cref{thm:SRW2}, but we omit this for simplicity.

\subsection{Organisation}

The rest of this paper is organised as follows. The proofs appear in the same order as the statements of the corresponding theorems in the introduction. The main results regarding the simple random walk are proved in \Cref{sec:SRW}. In \Cref{sec:variants}, we establish the variants of the SRW source localisation. \Cref{sec:conclude} contains some concluding remarks and open problems. We also include \Cref{app:cd,app:strongtrans}, concerning (A) estimates for the probability of a cut-edge in the two-sided simple random walk on $\ZZ^d$, and (B) strongly transient walks, including the equivalence of two different definitions used in the literature.

\section*{Acknowledgements}

We thank Roberto I. Oliveira for introducing us to the questions considered in this paper, during the Randomness and Learning on Networks (RandNET) program in 2024 at IMPA, Rio de Janeiro. We are grateful to IMPA for their hospitality, and also to Oliver Riordan and Jonathan Hermon for helpful discussions. RG is supported by a joint Clarendon Fund and Exeter College SKP scholarship. PK is supported by ERC Advanced Grant 883810. TP is supported by the Additional Funding Programme for Mathematical Sciences, delivered by EPSRC (EP/V521917/1) and the Heilbronn Institute for Mathematical Research.

\section{Proofs of the main results}
\label{sec:SRW}

In this section, we prove our main results regarding source localisation for the simple random walk. We begin by setting up the notation used throughout the paper.

\subsection{Notation}

Let $\gH$ be an infinite vertex-transitive strongly transient graph of finite degree. Let $(\overline{X}_n)_{n \in \ZZ}$ be the two-sided simple random walk on $\gH$, which can be defined as
\begin{equation*}
    \overline{X}_n = \begin{cases}
        X_n, \quad &n \ge 0,\\
        X'_{-n}, \quad &n < 0,
    \end{cases}
\end{equation*}
where $(X_n)_{n \ge 0}$ and $(X'_n)_{n \ge 0}$ are two independent copies of the simple random walk on $\gH$ started at the same (deterministic) vertex. For any interval $\mathcal{I}$ on the real line, we define $\overline{X} \mathcal{I} \coloneqq \{\overline{X}_n: n \in \mathcal{I} \cap \ZZ\}$ and write $G\mathcal{I}$ for the graph with $V(G\mathcal{I}) = \overline{X} \mathcal{I}$ and $E(G\mathcal{I}) = \{\{\overline{X}_n, \overline{X}_{n+1}\}: \{n,n+1\} \subseteq \mathcal{I}\}$.

We call $\{\overline{X}_n, \overline{X}_{n+1}\}$ a \emph{cut-edge} (of the two-sided random walk) if $\overline{X}(-\infty, n] \cap \overline{X}[n+1,\infty) = \emptyset$. We let $I_n$ denote the indicator of the event that $\{\overline{X}_n, \overline{X}_{n+1}\}$ is a cut-edge. We write $L$ for the set of cut-edges. We note that any cut-edge (of the two-sided random walk) is a cut-edge of the graph $G = G(-\infty, \infty)$; the reverse implication is not true in general. We define
\begin{equation}
\label{eqn:cH}
    \cgH \coloneqq \EE[I_0] = \PP(\overline{X}(-\infty, 0] \cap \overline{X}[1,\infty) = \emptyset).
\end{equation}
We also write $c(d) = c(\ZZ^d)$.
As $(I_n)_{n \in \ZZ}$ is a stationary process, the ergodic theorem implies
\begin{equation}
\label{eqn:ergodic}
    \lim_{n \to \infty} \frac{1}{n} \sum_{i=0}^{n-1} I_n = \EE[I_0] = \PP(I_0 = 1)
    \ \ \text{ almost surely.}
\end{equation}

For a finite connected graph $H$ with diameter $\diam(H)$, we let $\Diam(H)$ denote the set of diametrical paths in $H$ (paths of length $\diam(H)$ joining vertices at distance $\diam(H)$). We consider the following events on increases in diameter and the set $L$ of cut-edges.
\begin{align*}
    \cC^i_j &= \{ \diam(G[i,j]) < \diam(G[i,j+1]) \}, \\
    \cM^i_j &= \cM^i_j(\delta) =  \left\{ |E(G[i,j]) \cap L| \in \cgH (j-i) \cdot [1 - \delta, 1 + \delta] \right\}, \\
    \cN^i_j &= \cN^i_j(\eta) = \{ |E(P) \cap L| \ge (1 - \eta) \cgH (j-i), \forall P \in \Diam(G[i,j]) \}.
\end{align*}
Here $i \le j$ are integers and $\delta, \eta \in (0,1)$. For simpler notation, we also write $G_n = G[0,n]$, $\cC_n = \cC^0_n$, $\cM_n = \cM^0_n$, $\cN_n = \cN^0_n$, for $n \ge 0$. Intuitively, $\cC_n$ is the event that the diameter of $G_n$ increases in the next step, $\cM_n$ is the event that $G_n$ has (up to a constant close to $1$) as many cut-edges as expected, and $\cN_n$ is the event that every path in $\Diam(G_n)$ has many cut-edges.

\subsection{Localisation for strongly transient walks}
\label{sec:strongTransWalk}

In this section, we prove that positive constant probability localisation is possible for any strongly transient vertex-transitive graph. We recall from \Cref{sec:overview} that by reversibility it is equivalent to locate the final vertex $X_n$ rather than the source vertex $X_0$, and that to achieve this it suffices to establish a lower bound on $\liminf_n \PP(\cC_n)$, where $\cC_n$ is the event that step $n+1$ increased the diameter of the trace graph.

The proof is organised as follows. In \Cref{lem:Mn,lem:PCces}, we show concentration of the number of cut-edges and deduce that the Ces\`aro averages $\frac1n \sum_{i=0}^{n-1} \PP(\cC_i)$ are bounded away from $0$. Next, in \Cref{lem:Nn}, we show that with high probability every diametrical path of the trace has a large intersection with the set of cut-edges. This will be crucial in showing \Cref{lem:PCdiff}, which controls differences such as $|\PP(\cC_{2n})-\PP(\cC_n)|$. In \Cref{lem:cesaro} and \Cref{prop:PC}, we apply this control of linear scale variations to boost the Ces\`aro bounds to liminf bounds. We conclude by combining these ingredients in the \hyperref[prf:thmSRW1]{proof} of \Cref{thm:SRW}.

Our first lemma establishes the concentration of the number of cut-edges in $G_n$.

\begin{lemma}
\label{lem:Mn}
    For fixed $\delta \in (0, 1)$, we have $\lim_{n\to\infty} \PP(\cM_n(\delta)) = 1$.
\end{lemma}

\begin{proof}
    Recall from \eqref{eqn:ergodic} that almost surely $\lim_{n \to \infty} \frac{1}{n} \sum_{i=0}^{n-1} I_n = \EE[I_0] = \cgH$. Since almost sure convergence implies convergence in probability, we deduce
    \begin{equation*}
        \lim_{n \to \infty} \PP \Bigg( \Big|\frac{1}{n} \sum_{i=0}^{n-1} I_i - \cgH\Big| \le \delta \cgH\Bigg) = 1.
    \end{equation*}
    Observing that $\sum_{i=0}^{n-1} I_i = |E(G_n) \cap L|$ yields the desired result.
\end{proof}

Next we deduce the Ces\`aro bound for the events $\cC_n$ that the diameter of the trace increases.

\begin{lemma}
\label{lem:PCces}
    We have
    \begin{equation*}
        \liminf_{n\to\infty} \frac{1}{n} \sum_{i=0}^{n-1} \PP(\cC_i) \ge \cgH.
    \end{equation*}
\end{lemma}

\begin{proof}
    We begin with the simple but crucial observation that $\diam(G_{n+1}) \le \diam(G_n) + 1$, i.e., the diameter of the trace increases by at most one in a single step, so
    \begin{equation*} 
        \diam(G_n) \le \sum_{i=0}^{n-1} \ind_{\cC_i}, \text{ for all } n \in \ZZ_+.
    \end{equation*}
    Now, it follows from the definition of cut-edges that
    \begin{equation*}
        \diam(G_n) \ge d_{G_n}(X_0, X_n) \ge |E(G_n) \cap L|.
    \end{equation*}
    On the event $\cM_n(\delta)$, we have  $|E(G_n) \cap L| \ge \cgH (1-\delta) n$, so
    \begin{equation*}
        \sum_{i=0}^{n-1} \PP(\cC_i) \ge \EE[|E(G_n) \cap L|] \ge \cgH (1-\delta) n \cdot \PP(\cM_n).
    \end{equation*}
    Dividing both sides by $n$, taking $\liminf_{n \to \infty}$ and using \Cref{lem:Mn}, and subsequently taking $\lim_{\delta \to 0}$ yields the desired result.
\end{proof}

The next lemma controls the distribution of cut-edges in diametrical paths.

\begin{lemma}
\label{lem:Nn}
    For any fixed $\eta \in (0, 1)$, we have $\lim_{n\to\infty} \PP(\cN_n(\eta)) = 1$.
\end{lemma}

\begin{proof}
    Let $k = \lceil 16 \eta^{-1} \cgH^{-1} \rceil$. For each integer $n \ge 4 k/\eta$, we will define an event $A_n$ such that $A_n \subseteq \cN_n$, and $A_n$ occurs with high probability, which will imply the desired result.
    
    We begin by cutting the trace $X[0,n]$ into $k$ roughly equal sized parts $X[i_0,i_1]$, $X[i_1,i_2]$, $\ldots$, $X[i_{k-1},i_k]$, where $i_0 = 0$ and $i_k = n$. Here, the first $r$ parts have $m+1$ edges each, and the remaining $k-r$ parts have $m$ edges each, where $m = \lfloor n/k \rfloor$ and $r\in [0,k-1]$. Let $\delta = \eta/2$. Applying \Cref{lem:Mn}, we obtain that the number of cut-edges in the $i$th part lies in the interval $[(1-\eta/2) \cgH \cdot m, (1+\eta/2) \cgH \cdot (m+1)]$ with high probability. Let $A_n$ be the event that there are at least $(1 - \eta/2) \cgH \cdot m$ cut-edges in each part. Taking a union bound over $i \in [k]$, we conclude that the event $A_n$ occurs with high probability.

    It remains to show that $A_n \subseteq \cN_n$. Suppose that $A_n$ occurs. 
    Let $C_0$, $\ldots$, $C_\ell$ be the connected components of the graph obtained by deleting the set $L$ of cut-edges from $G[0,n]$; we list these components in a temporally ascending order. 
    Then each $C_i$ has at most $2m$ edges,
    since each $X[i_a, i_{a+1}]$ part has at least one cut-edge,
    so $C_i$ is contained in two consecutive parts.
            
    Let $u$ and $v$ be vertices with $d_{G_n}(u,v) = \diam(G_n)$. 
    Suppose that $u \in C_i$ and $v \in C_j$, where $i \le j$ without loss of generality. 
    We note that $i < j$, as
    \begin{equation*}
        d_{G_n}(u,v) = \diam(G_n) \ge |L \cap E(G_n)| \ge k (1-\eta/2) \cgH m > 2m \ge |E(C_i)|.
    \end{equation*}   
    Further, we claim that $i \le 2m$. For the proof, let $e$ denote the cut-edge between $C_i$ and $C_{i+1}$, and let $w$ be the endpoint of $e$ contained in $C_i$. Then we have
    \begin{align*}
        \diam(G_n) = d_{G_n}(u,v) &= d_{G_n}(u,w) + d_{G_n}(w,v)\\
        &= d_{C_i}(u,w) + d_{G_n}(w,v)\\
        &\le |E(C_i)| + d_{G_n}(w,v)\\
        &\le 2m + d_{G_n}(w,v).
    \end{align*}
    On the other hand, we have
    \begin{equation*}
        \diam(G_n) \ge d_{G_n}(X_0,v) = d_{G_n}(X_0,w) + d_{G_n}(w,v) \ge i + d_{G_n}(w,v).
    \end{equation*}
    Joining the two above inequalities yields the desired $i \le 2m$. Similarly, we have $j \ge \ell - 2m$. 
    
   Fixing any path $P$ of length $\diam(G_n)$ connecting $u$ and $v$ in $G_n$,
   we deduce that  $P$ contains all but $\ell - (j - i) \le 4m$ cut-edges. Thus,
    \begin{align*}
        |E(P) \cap L| &\ge k (1 - \eta/2) \cgH m - 4m\\
        &\ge \left(\frac{n}{k} - 1\right) (k (1 - \eta/2) \cgH - 4)\\
        &\ge \cgH n \left(1 - \frac{\eta}{2} - \frac{4}{k \cgH} - \frac{k}{n}\right)\\
        &\ge \cgH n \left(1 - \frac{\eta}{2} - \frac{\eta}{4} - \frac{\eta}{4}\right)\\
        &= \cgH (1-\eta) n. \qedhere
    \end{align*}
\end{proof}

Our next lemma shows that the event $\cC^{n-i}_n$ that step $n+1$ increases the diameter of $G[n-i,n]$ is not very sensitive to the choice of $i$.

\begin{lemma}
\label{lem:PCdiff}
    We have $\sup_{i \in [n/3,n]} \PP(\cC^0_n \symdiff \cC^{n-i}_n) \xrightarrow[n\to\infty]{} 0$, where $\symdiff$ is the symmetric difference.
\end{lemma}

\begin{proof}
    Fix $\delta$, $\eta \in (0, 0.1)$ and consider the following event that controls cut-edges
    in (a) diametrical paths of $G_n$, (b) diametrical paths of $G^{n-i}_n$, and (c) the full trace $G_n$:
    \begin{equation*}
        \cA_{n,i} = \cA_{n,i}(\delta,\eta) = \cN^0_{n}(\eta) \cap \cN^{n-i}_{n}(\eta) \cap \cM^0_{n}(\delta).
    \end{equation*}
    By \Cref{lem:Mn,lem:Nn}, using $X[n-i,n] \stackrel{d}{=} X[0,i]$, we have
    \begin{equation*}
        \sup_{i \in [n/3,n]} \PP(\cA_{n,i}^\complement) \le \PP(\cN_n^\complement) + \sup_{i \in [n/3,n]} \PP(\cN_i^\complement) + \PP(\cM_n^\complement) \xrightarrow[n\to\infty]{} 0.
    \end{equation*}
    We will prove that $\cC^0_n \symdiff \cC^{n-i}_n \subseteq \cA_{n,i}^\complement$, which will imply the lemma.
    
    Fix some $n \in \NN$ and $i \in [n/3,n]$. Suppose $P_1$ and $P_2$ are any paths in $\Diam(G[0,n])$ and $\Diam(G[n-i,n])$, respectively. Recalling that $L$ denotes the set of cut-edges, we define
    \begin{equation*}
        L_1 = E(P_1) \cap L, \quad L_2 = E(P_2) \cap L.
    \end{equation*}
    
    \begin{claim}
    \label{cla:dis}
        Whenever $\cA_{n,i}$ occurs, we have $L_1 \cap L_2 \neq \emptyset$.
    \end{claim}
    \begin{proof}[Proof of \Cref{cla:dis}] 
    It follows from the definitions of $\cN^0_n$ and $\cN^{n-i}_n$ that
    \begin{equation*}
        |L_1| > 0.9 \cgH n, \text{ and } |L_2| > 0.9 \cgH i.
    \end{equation*}
   Suppose for contradiction that $L_1$ and $L_2$ are disjoint. Then
    \begin{equation*}
        |G[0,n] \cap L| \ge |L_1| + |L_2| > 0.9 \cgH n + 0.9 \cgH i \ge 1.1 \cgH n.
    \end{equation*}
    However, this contradicts the definition of $\cM^0_n$, so the claim follows.
    \end{proof}

    The following claim will complete the proof of the lemma, as it implies
    \begin{equation*}
        \PP(\cC^0_n \symdiff \cC^{n-i}_n) \le \PP\big( (\cC^0_n \cap \cA_{n,i}) \symdiff(\cC^{n-i}_n \cap \cA_{n,i})\big) + \PP(\cA_{n,i}^\complement) = \PP(\cA_{n,i}^\complement). \phantom{\qedhere}
    \end{equation*}

    \begin{claim}
    \label{cla:int}
        We have $\cC^0_n \cap \cA_{n,i} = \cC^{n-i}_n \cap \cA_{n,i}$.
    \end{claim}
    \begin{proof}[Proof of \Cref{cla:int}]
    Suppose that $\cA_{n,i}$ occurs. We want to show that $\cC^0_n$ occurs if and only if $\cC^{n-i}_n$ occurs. 
    Suppose that $\cC^0_n$ occurs, i.e.~$\diam(G_{n+1}) > \diam(G_n)$. Fix any path $P \in \Diam(G_{n+1})$. Then $X_{n+1}$ must be one of the endpoints of $P$. Let $w$ be the other endpoint of $P$. Further, let $P'$ be the path obtained by deleting $X_{n+1}$ from $P$. Then, $P' \in \Diam(G_n)$. 
    
    Now let $Q$ be any path in $\Diam(G[n-i,n])$. Then it follows from \Cref{cla:dis} that $P'$ and $Q$ have a common cut-edge, say $\{u,v\}$. Without loss of generality, we may assume that $u$ appears before $v$ in the temporal order. Let $x$ and $y$ be the endpoints of $Q$ such that $x$ appears before $y$ in the temporal order. Since $\{u, v\}$ is a cut-edge, we have
    \begin{equation} \label{eq:diameters}
        d_{G[n-i,n]}(v, y) = d_{G_n}(v, y) \quad \text{ and } \quad d_{G[n-i,n]}(v, X_n) = d_{G_n}(v, X_n).
    \end{equation}
    Observe that
    \begin{align*}
        d_{G_n}(w, v) + d_{G_n}(v, y) &= d_{G_n}(w, y) \\ &\le \diam(G_n) \\ &= d_{G_n}(w, X_n) = d_{G_n}(w, v) + d_{G_n}(v, X_n),
    \end{align*}
    and so by \eqref{eq:diameters} we obtain $d_{G[n-i,n]}(v,y) \le d_{G[n-i,n]}(v,X_n)$.
    Further,
    \begin{align*}
        \diam(G[n-i, n+1]) &\ge d_{G[n-i, n+1]}(x, X_{n+1})\\
        &= 1 + d_{G[n-i,n]}(x, X_n)\\
        &> d_{G[n-i,n]}(x, X_n)\\
        &= d_{G[n-i,n]}(x, v) + d_{G[n-i,n]}(v, X_n)\\
        &\ge d_{G[n-i,n]}(x, v) + d_{G[n-i,n]}(v, y)\\
        &= d_{G[n-i,n]}(x, y)\\
        &= \diam(G[n-i,n]).
    \end{align*}
    Therefore, we conclude that the event $\cC^{n-i}_n$ occurs. 
    A similar argument gives the other direction,
    thus proving the claim, and so the lemma.
    \end{proof}
\end{proof}

Next we require a simple lemma from real analysis that will be used in \Cref{prop:PC} to boost Ces\`aro control to liminf control.

\begin{lemma}[Ces\`aro means]
\label{lem:cesaro}
    Suppose $(c_n)$ is a bounded sequence satisfying
    \begin{equation}
        \sup_{i \in [n/2,n]}|c_n - c_i| \xrightarrow[]{n\to\infty} 0. \label{eq:assump_ces}
    \end{equation}
    Then
    \begin{equation*}
        \liminf_{n\to\infty} c_n \ge \liminf_{n\to\infty}{\frac{1}{n} \sum_{i=1}^n c_i}.
    \end{equation*}
\end{lemma}

\begin{proof}
    By shifting all the elements, without loss of generality, we may assume that $0 \le c_n \le B$ for all $n$ for some $B \ge 0$. Further, let
    \begin{equation*}
        \varepsilon_n = \sup_{m \ge n} \sup_{i \in [m/2,m]} |c_m - c_i|.
    \end{equation*}
 We note that $(\varepsilon_n)$ is non-increasing and 
$\varepsilon_n \xrightarrow[]{} 0$ as $n \to \infty$ by \eqref{eq:assump_ces}.
    
    Fix an arbitrary $k \ge 1$. Note that for every $n \in [k, 2k]$, we have
    \begin{equation*}
        c_n \le c_k + \varepsilon_n \le c_k + \varepsilon_k.
    \end{equation*}
    By inductive use of the above inequality, as $(\varepsilon_n)$ is non-increasing, we get
    \begin{equation}\label{eq:cnck}
        c_n \le c_{k} + j\varepsilon_{k},
    \end{equation}
    for all $j \ge 1$ and $n \in [k 2^{j-1}, k 2^j]$. As $\varepsilon_k \ge 0$, it follows that \eqref{eq:cnck} holds for all $n \in [k, k 2^j]$.
    
    Let $s_n = \frac{1}{n}\sum_{i=1}^n c_i$. As $(c_n)$ is bounded and non-negative, by \eqref{eq:cnck} we have
    \begin{align*}
        s_{k2^j} \le \frac{1}{k2^j}\left(kB + \sum_{i=k+1}^{k2^j} (c_k + j \varepsilon_k)\right) \le \frac{B}{2^j} + c_k + j \varepsilon_k.
    \end{align*}
    Taking $\liminf_{k\to\infty}$ and subsequently $\liminf_{j \to \infty}$ yields the statement.
\end{proof}

\begin{proposition}
\label{prop:PC}
    We have
       $ \liminf_{n \to \infty} \PP(\cC_n) \ge \cgH$.
\end{proposition}

\begin{proof}
The statement follows from \Cref{lem:cesaro,lem:PCces} with $c_n = \PP(\cC_n)$, noting that
    \begin{equation*}
        \PP(\cC^0_n \symdiff \cC^{n-i}_n) \ge |\PP(\cC^0_n) - \PP(\cC^{n-i}_n)| = |c_n - c_i|,
    \end{equation*}
    and so by \Cref{lem:PCdiff} the assumption \eqref{eq:assump_ces} holds. 
 \end{proof}

We now have all the ingredients to prove the main result of this section, \Cref{thm:SRW}.

\begin{proof}[Proof of \Cref{thm:SRW}]
\label{prf:thmSRW1}
    Let $n \in \ZZ_+$. Suppose that the event $\cC_{n-1}$ occurs. 
    Then each pair $(\bfu, \bfv)$ of vertices in $G_n$ with $d_{G_n}(\bfu, \bfv) = \diam(G_n)$ contains $X_n$.
    As $\Psi(G)$ selects uniformly at random one endpoint of some diameter of $G$, 
    by reversibility we have
    \begin{equation*}
        \PP(\Psi(G_n) = X_0) = \PP(\Psi(G_n) = X_n) \ge \frac{1}{2} \PP(\cC_{n-1}).
    \end{equation*}
    Taking $\liminf_{n \to \infty}$, the result follows using \Cref{prop:PC} and $\cgH > 0$ by \Cref{lem:intprob}.
\end{proof}

\subsection{Optimal bounds in large dimensions}

In this section, we prove \Cref{thm:SRWoptimal} concerning optimal localisation bounds for $\gH = \ZZ^d$ as $d\to\infty$. The estimator $\Psi$ used in the proof of \Cref{thm:SRW} has localisation probability $c(d)/2$, which is $1/2 - \Theta(1/d)$ by \Cref{lem:cd}. Here we will prove a matching upper bound on the localisation probability of any estimator. In \Cref{sec:overview}, we indicated how one can achieve this via coupling arguments and sketched a simple upper bound of $1/2$ using the involution of reversing walks. Here we will combine this with another coupling to improve this to $1/2 - \Omega(1/d)$.

\begin{figure}[h!]
    \centering
    \includegraphics[width=0.7\linewidth]{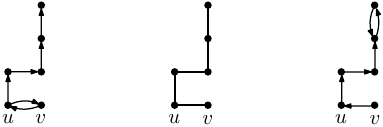}
    \caption{Two walks producing the same trace. The left and right walks start at $u$ and $v$, respectively.}
    \label{fig_couplingfinite}
\end{figure}

Observe that if the beginning of the walk looks as in the left picture of \Cref{fig_couplingfinite}, then one cannot distinguish between vertices $u$ and $v$: there is a corresponding walk in the right picture of \Cref{fig_couplingfinite}, which starts at a different vertex but produces the same trace. The probability of the walk starting as in the figure is $\Theta(1/d)$. On this event, which has probability $p = \Theta(1/d)$, we bound the localisation probability by $1/4$ rather than $1/2$ (which comes from the coupling of reversing walks), so our overall bound is $\tfrac14p + \tfrac12(1-p) = \tfrac12 - \Omega(1/d)$. The argument for infinite traces is simpler, as we do not need to consider time reversal; see the \hyperref[prf:upperboundinfinite]{proof} of \Cref{thm:infinite}. 

\begin{proof}[Proof of the upper bound in \Cref{thm:SRWoptimal}]
\label{prf:ub}
    Let $E \coloneqq \{\bfe_1, \dots, \bfe_d\}$ be the set of standard basis vectors for $\ZZ^d$, and let $E_\pm \coloneqq \{\pm \bfe: \bfe \in E\}$. We shall henceforth assume $n \ge 3$. Let $\Omega \coloneqq E_\pm^n$ be the space of sequences of length $n$ with entries in $E_\pm$. We define maps $f, g, h: \Omega \to \Omega$ by
    \begin{align*}
        f((s_1,\ldots,s_n)) &= (-s_n, \ldots, -s_1), \\
        g((s_1,\ldots,s_n)) &= \begin{cases}
            (s_2, \ldots, s_n, -s_n), &\text{ if } s_1 = -s_2, \\
            (s_2, \ldots, s_n, -s_2), &\text{ if } s_1 \ne -s_2 \text{ and } s_1 = -s_n,\\
            (s_2, \ldots, s_n, s_1), &\text{ otherwise},
        \end{cases}\\
        h((s_1,\ldots,s_n)) &= (s_n, s_1, \ldots, s_{n-1}).
    \end{align*}
    Let $S = (S_1, \dots, S_n)$ be a sequence of iid random variables with uniform distribution on $E_\pm$. Observe that $f \circ f$ and $h \circ g \circ h \circ g$ are both equal to the identity map on $\Omega$. Consequently, all of $f$, $g$, and $f\circ g$ are bijections on $\Omega$. Note that the law of $S$ is uniform on $\Omega$. Therefore, the laws of $f(S)$, $g(S)$, and $f(g(S))$ are all uniform on $\Omega$. Let $X$, $X^f$, $X^g$, and $X^{fg}$ be $n$-step simple random walks on $\ZZ^d$ started at the origin with steps given by $S$, $f(S)$, $g(S)$ and $f(g(S))$, respectively. Further, let $G_n$ be the trace graph of the walk $X$ embedded inside $\ZZ^d$, and similarly define $G_n^f$, $G_n^g$, and $G_n^{fg}$.

    For $\ell \in \NN$, let us define the box $B_\ell \coloneqq [-\ell, \ell]^d \cap \ZZ^d$ and consider the random variable
    \begin{align*}
        W &\coloneqq \sum_{\bfv \in B_\ell} \mathbf{1}\{\Phi(G_n + \bfv) = \bfv\} + \sum_{\bfv \in B_\ell} \mathbf{1}\{\Phi(G^f_n + \bfv + X_n) = \bfv + X_n\}\\
        &+ \sum_{\bfv \in B_\ell} \mathbf{1}\{\Phi(G_n^g + \bfv + X_1) = \bfv + X_1\} + \sum_{\bfv \in B_\ell} \mathbf{1}\{\Phi(G^{fg}_n + \bfv + X_1 + X_n^g) = \bfv + X_1 + X_n^g\}.
    \end{align*}
    Similarly to the simplified sketch in \Cref{sec:overview}, we will bound $ \EE[W]$ below in \Cref{cla:Wlower} in terms of the localisation probability and above in \Cref{cla:Wupper} by considering the events that certain indicators in the definition of $W$ cannot occur simultaneously, as they would require the estimator to choose two different vertices for the same trace.
    
    \begin{claim}
    \label{cla:Wlower}
        We have $\EE[W] \ge (|B_\ell| + |B_{\ell-n}| + |B_{\ell-1}| + |B_{\ell-n-1}|) \inf_{\bfv \in \ZZ^d} \PP(\Phi(G_n + \bfv) = \bfv)$.
    \end{claim}
    
    \begin{proof}[Proof of \Cref{cla:Wlower}]
        Taking expectation on both sides, we obtain
        \begin{align}
            \label{eqn:expWf}
            \EE[W] &= \sum_{\bfv \in B_\ell} \PP(\Phi(G_n + \bfv) = \bfv) + \sum_{\bfv \in B_\ell} \PP(\Phi(G^f_n + \bfv + X_n) = \bfv + X_n) \\
            &+ \sum_{\bfv \in B_\ell} \PP(\Phi(G_n^g + \bfv + X_1) = \bfv + X_1) + \sum_{\bfv \in B_\ell} \PP(\Phi(G^{fg}_n + \bfv + X_1 + X_n^g) = \bfv + X_1 + X_n^g). \nonumber
        \end{align}
        The second sum on the right-hand side above can be rewritten as
        \begin{equation*}
            \sum_{\bfv \in B_\ell} \sum_{\bfw \in B_n} \PP(X_n = \bfw)\, \PP(\Phi(G_n^f + \bfv + \bfw) = \bfv + \bfw \mid X_n = \bfw),
        \end{equation*}
        which can further be rewritten as
        \begin{equation*}
            \sum_{\bfv' \in B_{\ell+n}} \sum_{\bfw \in B_n} \mathbf{1}\{\bfv'-\bfw \in B_\ell\}\, \PP(X_n = \bfw)\, \PP(\Phi(G_n^f + \bfv') = \bfv' \mid X_n = \bfw).
        \end{equation*}
        Now observing that $\bfv' - \bfw \in B_\ell$ for all $\bfv' \in B_{\ell - n}$ and $w \in B_n$, the above sum is at least
        \begin{equation*}
            \sum_{\bfv' \in B_{\ell-n}} \sum_{\bfw \in B_n} \PP(X_n = \bfw)\, \PP(\Phi(G_n^f + \bfv') = \bfv' \mid X_n = \bfw) = \sum_{\bfv' \in B_{\ell-n}} \PP(\Phi(G_n^f + \bfv') = \bfv').
        \end{equation*}
        By a similar argument as above, we conclude that the third and fourth sums on the right-hand side of \eqref{eqn:expWf} are bounded below by
        \begin{equation*}
            \sum_{\bfv' \in B_{\ell-1}} \PP(\Phi(G_n^g + \bfv') = \bfv'), \quad \text{ and } \quad \sum_{\bfv' \in B_{\ell-n-1}} \PP(\Phi(G_n^{fg} + \bfv') = \bfv'),
        \end{equation*}
        respectively. Using the above estimate in \eqref{eqn:expWf}, we obtain
        \begin{align*}
            \EE[W] \ge &\sum_{\bfv \in B_\ell} \PP(\Phi(G_n + \bfv) = \bfv) + \sum_{\bfv' \in B_{\ell-n}} \PP(\Phi(G_n^f + \bfv') = \bfv')\\
            &+ \sum_{\bfv' \in B_{\ell-1}} \PP(\Phi(G_n^g + \bfv') = \bfv') + \sum_{\bfv' \in B_{\ell-n-1}} \PP(\Phi(G_n^{fg} + \bfv') = \bfv'),
        \end{align*}
        which in turn implies the claim.
    \end{proof}
  
    \begin{claim}
    \label{cla:Wupper}
        There exists a sequence $(\varepsilon_n)_{n \ge 6}$ of nonnegative reals converging to $0$ such that $\EE[W] \le (2 - 1/(2d) + \varepsilon_n) |B_\ell|$ for all $\ell, n \in \NN$ with $n \ge 6$.
    \end{claim}
    
    \begin{proof}[Proof of \Cref{cla:Wupper}]
        We consider the events    
        \begin{align*}
            A &= \{S_1=-S_2\}, \\
            B &= \{ \| S_3+S_4+\ldots +S_{n-1}\|_1 > 5\}, \\
            C &= A \cap B, \\
            D &= A^\complement \cap B.
        \end{align*}
        Note that the walks $X^f$ and $X^{fg}$ are the time reversals of the walks $X$ and $X^g$, respectively. Therefore, we have
        \begin{equation*}
            G_n = G_n^f + X_n, \quad \text{ and } \quad G_n^g = G_n^{fg} + X_n^g.
        \end{equation*}
        Observe that whenever $D$ occurs, we have $X_n \ne 0$ and $X_n^g \ne 0$. Thus, whenever $D$ occurs, for $\bfv \in \ZZ^d$, we have
        \begin{equation*}
            \mathbf{1}\{\Phi(G_n + \bfv) = \bfv\} + \mathbf{1}\{\Phi(G^f_n + \bfv + X_n) = \bfv + X_n\} \le 1,
        \end{equation*}
        and
        \begin{equation*}
            \mathbf{1}\{\Phi(G_n^g + \bfv + X_1) = \bfv + X_1\} + \mathbf{1}\{\Phi(G^{fg}_n + \bfv + X_1 + X_n^g) = \bfv + X_1 + X_n^g\} \le 1.
        \end{equation*}
        Adding the two inequalities above and then summing over $\bfv \in B_\ell$ implies that $W \le 2|B_\ell|$ whenever $D$ occurs. Similarly, observe that whenever $C$ occurs, we have
        \begin{equation*}
            G_n = G_n^f + X_n = G_n^g + X_1 = G_n^{fg} + X_1 + X_n^g,
        \end{equation*}
        and the points $0$, $X_n$, $X_1$, and $X_1 + X_n^g$ are distinct. Thus, whenever $C$ occurs, we have
        \begin{align*}
            \mathbf{1}\{\Phi(G_n &+ \bfv) = \bfv\} + \mathbf{1}\{\Phi(G^f_n + \bfv + X_n) = \bfv + X_n\}\\
            &+ \mathbf{1}\{\Phi(G_n^g + \bfv + X_1) = \bfv + X_1\} + \mathbf{1}\{\Phi(G^{fg}_n + \bfv + X_1 + X_n^g) = \bfv + X_1 + X_n^g\} \le 1,
        \end{align*}
        for $\bfv \in \ZZ^d$. Summing the above over $\bfv \in B_\ell$ implies that $W \le |B_\ell|$ whenever $C$ occurs. Finally, it follows from the definition of $W$ that $W \le 4 |B_\ell|$. Hence, we obtain
        \begin{align*}
            \EE[W] &= \PP(B^\complement) \EE[W|B^\complement] + \PP(C) \EE[W|C] + \PP(D) \EE[W|D]\\
            &\le \PP(B^\complement) \cdot 4 \cdot |B_\ell| + \PP(A) \cdot |B_\ell| + \PP(A^\complement) \cdot 2 \cdot |B_\ell|\\
            &= (2 - \PP(A) + 4 \PP(B^\complement)) |B_\ell|.
        \end{align*}
        The claim follows as  $\PP(A) = 1/(2d)$ and $\PP(B^\complement) \to 0$ as $n \to \infty$.
    \end{proof}
    
    Combining \Cref{cla:Wlower} and \Cref{cla:Wupper}, we obtain
    \begin{equation*}
        (|B_\ell| + |B_{\ell-n}| + |B_{\ell-1}| + |B_{\ell-n-1}|) \inf_{\bfv \in \ZZ^d} \PP(\Phi(G_n + \bfv) = \bfv) \le (2 - 1/(2d) + \varepsilon_n) |B_\ell|.
    \end{equation*}
    Observe that $|B_\ell| = \Theta(\ell^d)$ and $|B_{\ell-n}|$, $|B_{\ell-1}|$, $|B_{\ell-n-1}|$ are all $|B_\ell| - O(n \cdot \ell^{d-1})$ as $\ell \rightarrow \infty$. Therefore, dividing both sides of the above inequality by $|B_\ell| + |B_{\ell-n}| + |B_{\ell-1}| + |B_{\ell-n-1}|$, and taking $\lim_{\ell \to \infty}$, we obtain
    \begin{equation*}
        \inf_{\bfv \in \ZZ^d} \PP(\Phi(G_n + \bfv) = \bfv) \le 1/2 - 1/(8d) + \varepsilon_n/4.
    \end{equation*}
    Taking $\limsup_{n \to \infty}$ completes the proof.
\end{proof}

\subsection{No localisation for recurrent walks}
\label{sec:recWalk}

In this section, we prove that positive constant probability localisation is not possible for recurrent graphs. The main idea of the proof is the same as in the \hyperref[prf:ub]{proof} of the upper bound in \Cref{thm:SRW2}, namely, to couple walks started at distinct vertices so that their traces are exactly the same. However, in this case, not only can we find a high probability coupling with the desired property, we can also find it for a growing number of starting points. Since the estimator cannot distinguish between the starting points, its success rate tends to $0$ as the number of walks increases. 

For the proof, we will need a lemma saying that simple random walks on a recurrent graph eventually collide, subject to a
mild parity condition that is necessary in case the graph is bipartite. The proof follows from Proposition~2.1 and Remark~1.1 in \cite{BPS12}. Note also that the conclusion of the lemma may not hold if $\gH$ is not vertex-transitive; see \cite{KP04}.

\begin{lemma}
\label{lem:finitecollision}
    Suppose $\gH$ is an infinite, connected, vertex-transitive, recurrent graph of finite degree. Further, let $a$, $b \in \gH$ be such that the transition probability $p_{2j}(a, b)$ is positive for some $j \ge 0$, and consider two independent simple random walks $X_t$, $Y_t$ on $\gH$. Then, we have
    \begin{equation*}
        \PP(\exists t > 0 \colon X_t = Y_t \mid X_0 = a, Y_0 = b) = 1.
    \end{equation*}
\end{lemma}

\begin{proof}[Proof of \Cref{thm:SRWrecurrent}]
    As indicated in \Cref{sec:overview}, we consider a large number $k$ of random walks starting from distinct vertices, which we aim to couple so that they have identical traces, and so cannot be distinguished by any estimator.
    
    Fix some arbitrary $v \in \gH$ and let $v_1, \dots, v_k$ be vertices of $\gH$ such that for each $i \in [k]$, there exists $j > 0$ such that $p_{2j}(v,v_i) > 0$. Let $X, X^1$, $\ldots$, $X^k$ be simple random walks on $\gH$ such that $X_0 = v$ and $X^i_0 = v_i$ for $i \in [k]$. We couple the walks in the following standard manner: each walk $X^i$ performs steps independently of other walks until it meets $X$ (at the same vertex at the same time $t_0$). From that point onwards, $X^i$ follows the steps of $X$ so that for every $t \ge t_0$ we have $X^i_{t} = X_{t}$.

    By \Cref{lem:finitecollision}, almost surely all the walks $X^i$ eventually couple with $X$. In particular, there is a constant $t_1 < \infty$ such that $\PP(A) \ge 1 - 1/k$, where
    \begin{equation*}
        A = \bigcap_{i \in [k]} \{X^i_{t_1} = X_{t_1}\}.
    \end{equation*}
    Further, as $X$ is recurrent, there is a constant $t_2 < \infty$ such that $\PP(B) \ge 1 - 1/k$, where
    \begin{equation*}
        B = \{\{(X_{t_1+1}, X_{t_1+2}), (X_{t_1+2}, X_{t_1+3}), \ldots, (X_{t_1+t_2-1}, X_{t_1+t_2})\} \supseteq E(B_{t_1}(X_{t_1}))\}
    \end{equation*}  
    and $E(B_r(x))$ is the set of edges of the ball centered at $x$ of radius $r$.
    
    Let $T_n(Y)$ be the trace graph of the first $n$ steps of the walk $Y$. For $n \ge t_1+t_2$ we have
    \begin{equation}
    \label{eq:equaltraces}
        \PP(\forall i \in [k] \colon \ T_n(X^i) = T_n(X)) \ge \PP(A \cap B) \ge  1 - 2/k.    
    \end{equation}
    Indeed, the event $A$ ensures that the traces of walks $X^i$ after time $t_1$ coincide with the trace of $X$, while the event $B$ ensures that $X$ traverses the traces of $X^i$'s (and also its own) created before time $t_1$. As a consequence, all their total traces created until time $n$ coincide.
    
    By \eqref{eq:equaltraces}, for $n \ge t_1+t_2$, we have
    \begin{align*}
        \sum_{i=1}^k \PP_{v_i}(\Phi(G_n^\gH) = v_i) = \sum_{i=1}^k \PP(\Phi(T_n(X^i)) = v_i)
        &\le 2 + \sum_{i=1}^k \PP(A \cap B \cap \{\Phi(T_n(X^i)) = v_i\})\\
        &= 2 + \sum_{i=1}^k \PP(A \cap B \cap \{\Phi(T_n(X)) = v_i\}) \le 3,
    \end{align*}
    where we used the fact that on $A\cap B$ one has $T_n(X^i) = T_n(X)$, and the fact that all events $\{\Phi(T_n(X)) = v_i\}$ are disjoint due to the $v_i$ being distinct. This implies
    \begin{align*}
        \limsup_{n\to\infty} \inf_{w \in \gH} \PP_w(\Phi(G_n^\gH) = w) \le \limsup_{n\to\infty} \frac{1}{k} \sum_{i=1}^k \PP_{v_i}(\Phi(G_n^\gH)=v_i) \le \frac{3}{k}.
    \end{align*}
    Finally, taking $k \to \infty$ yields the statement.
\end{proof}

\section{Variants of localisation for the SRW}
\label{sec:variants}

In this section, we will consider different variants of the source localisation problem for the simple random walk. For simplicity, we fix $\gH=\ZZ^d$, although our arguments apply to any strongly transient host graph $\gH$. We remind the reader that $c(d) = c(\ZZ^d) > 0$ for $d \ge 5$ by \Cref{lem:cd}. \Cref{sec:const} considers the high accuracy problem where we aim to output a set of vertices that contains the source with probability close to $1$. The remaining sections consider the usual single vertex estimators but with variations on the input graph: the infinite trace in \Cref{sec:infinite}, the trace of exactly $n$ vertices in \Cref{sec:range}, and the vertex trace in \Cref{sec:vertex}.
 
\subsection{High accuracy source localisation}
\label{sec:const}

In this section, we consider the following estimator $\Lambda_k$ for locating the source in a set of size $k$. Given an input graph $G$ and $k \ge 2$ even, we define $\Lambda_k(G)$ by selecting a uniformly random pair $(u,v)$ of diametric vertices in $G$ and outputting the $k/2$ closest vertices to $u$ and the $k/2$ closest vertices to $v$.

To estimate the success probability of $\Lambda_k$, we will need to control the number of cut-edges in segments of the walk that are considerably smaller than those considered in \Cref{sec:strongTransWalk}.

For integers $i \le j$, let $L[i,j] = L \cap X[i,j]$ denote the set of cut-edges within the walk $X[i,j]$. 

We begin with an estimate on the variance of $|L_k|$, where $L_k \coloneqq L[0,k]$.

\begin{lemma}
\label{lem:var}
    For every $d \ge 5$, we have $\mathrm{Var}(|L_k|) = O(k^{3/2})$ as $k \xrightarrow{} \infty$.
\end{lemma}

\begin{proof}
    For $n \in \ZZ$, let $e_n$ be the edge $\{X_{n-1}, X_n\}$. Writing $|L_k|$ as a sum of indicators, we obtain
    \begin{equation}
    \label{eqn:var}
        \text{Var}(|L_k|) = \sum_{i=1}^k c(d) (1-c(d)) + 2 \sum_{1 \le i < j \le k} (\PP(\{e_i,e_j\} \subseteq L) - c(d)^2).
    \end{equation}
    We will now estimate the second summand on the right-hand side. For $1 \le i < j \le k$, suppose that $e_i$ and $e_j$ are both cut-edges. Cut the two-sided SRW trace into three segments, the left segment $\overline{X}(-\infty,i-1]$ before the edge $e_i$, the middle segment $\overline{X}[i,j-1]$ between the edges $e_i$ and $e_j$, and the right segment $\overline{X}[j,\infty)$ after the edge $e_j$. Then the left segment does not intersect the left half of the middle segment. Similarly, the right segment does not intersect the right half of the middle segment. Therefore, we have
    \begin{align}
    \label{eqn:summand}
        \PP(\{e_i,e_j\} \subseteq L) &\le \PP\left(\overline{X}(-\infty,i) \cap \overline{X}[i,\tfrac{i+j-1}{2}) = \overline{X}(\tfrac{i+j-1}{2}, j) \cap \overline{X}[j,\infty) = \emptyset)\right) \nonumber \\
        &= \PP\left(\overline{X}(-\infty,0] \cap \overline{X}[1, \lceil (j-i-1)/2 \rceil] = \emptyset\right)^2.
    \end{align}
    To see the equality above, we construct the two-sided random walk via two independent one-sided random walks starting in the middle of the middle segment.
    
    The probability on the right-hand side of \eqref{eqn:summand} is bounded above by
    \begin{align*}
        c(d) + \PP\left(\overline{X}(-\infty,0] \cap \overline{X}[\lceil (j-i-1)/2 \rceil + 1, \infty) \ne \emptyset\right).
    \end{align*}
    As in \Cref{lem:pd} and \Cref{lem:lclt}, the probability above can be in turn bounded by
    \begin{equation*}
        \sum_{n = 1}^\infty n \PP(X_{n + \lceil (j-i-1)/2 \rceil} = 0) = \sum_{n=1}^{\infty} \frac{O(1)}{(n+(j-i-1)/2)^{d/2-1}} = \frac{O(1)}{(j-i)^{d/2-2}} = \frac{O(1)}{(j-i)^{1/2}}.
    \end{equation*}
    Substituting the above estimates in \eqref{eqn:summand}, and then in \eqref{eqn:var}, proves the lemma.
\end{proof}

Next we will define good events that occur with reasonably large probability on which our estimator $\Lambda_k(G_n)$ contains the source vertex. Fix $m \in \ZZ^+$ with $m \ge 8/c(d)$, let $\lambda = c(d)/8$ and
\begin{equation*}
    a_0 = m, \quad a_1 = \left\lceil \lambda m \right\rceil, 
    \quad a_n = \left\lceil \lambda b_{n-2} \right\rceil, \text{ for all } n \ge 2.
\end{equation*}
Let $(b_n)_{n \ge -1}$ satisfy $b_{-1} = 0$ and $b_n = \sum_{m = 0}^{n} a_m$ for $n\ge 0$. We define
\begin{equation}
\label{eq:eventAm}
    \cA_m \coloneqq \bigcap_{n \ge 0} \left\{L[b_{n-1}, b_n] \ge \frac{c(d)}{2} a_n\right\}.
\end{equation}
Thus $\cA_m$ is the event that at least a $c(d)/2$ proportion of edges are cut-edges in each of the segments $G[b_{n-1},b_n]$. The choice of the sequences $(a_n)$ and $(b_n)$ was governed by the need to satisfy the following inequalities.

\begin{lemma}
\label{lem:cond}
    We have $a_n + a_{n+1} \le \frac{c(d)}{2} b_{n-1}$ and $a_n \ge \lambda \left(1 + \frac{\lambda}{1+\lambda}\right)^{n-2} m$ for all $n \ge 2$.
\end{lemma}

\begin{proof}
    As $c(d)m \ge 8$, we have
    \begin{equation*}
        a_1 + a_2 \le 2\lambda m + 2 \le \frac{c(d)}{2}m = \frac{c(d)}{2} a_0.
    \end{equation*}
    Furthermore, for any $n \ge 2$, we have
    \begin{equation*}
        a_n + a_{n+1} \le 2\lambda b_{n-2} + \lambda a_{n-1} + 2 \le \frac{c(d)}{2} b_{n-1}.
    \end{equation*}
    Now we show the second inequality by induction. The base case $n = 2$ is direct from the definition. Now assuming the inequality for $n \le \ell$, we have
    \begin{align*}
        a_{\ell+1} &\ge \lambda \left[m + \lambda m + \sum_{i=2}^{\ell-1} a_i \right]\\
        &\ge \lambda \left[m + \lambda m + \sum_{i=2}^{\ell-1} \lambda \left(1 + \frac{\lambda}{1+\lambda}\right)^{i-2} m \right]\\
        &= \lambda (1+\lambda) \left(1 + \frac{\lambda}{1+\lambda}\right)^{\ell-2} m\\
        &\ge \lambda \left(1 + \frac{\lambda}{1+\lambda}\right)^{\ell-1} m. \qedhere
    \end{align*}
\end{proof}

In the following lemma, we show that whenever $\cA_m$ occurs, the source is among the $O(m)$ closest vertices to one (or both) of the arbitrarily chosen two most distant points in $G_n$.

\begin{lemma}
\label{lem:agood}
    Let $m \ge 8/c(d)$. Suppose $\cA_m$ occurs. Let $P$ be any path in $\Diam(G_n)$. Then, the source is among the $\lfloor 2m(2/c(d)+1) \rfloor$ closest vertices to one (or both) of the endpoints of $P$.
\end{lemma}

\begin{proof}
    Suppose that $P$ has endpoints $u$ and $v$. Recall that $L$ denotes the set of cut-edges.
    Let $C_1, \dots, C_{\ell}$ be the connected components of $G_n \setminus L$ arranged in a temporally ascending order.

    \begin{claim}
    \label{cla:uv}
        At least one of $u$ and $v$ appears in the first segment $G[0, b_0]$.
    \end{claim}
    
    \begin{proof}[Proof of \Cref{cla:uv}]
        Suppose not. Further, suppose for contradiction that $u$ and $v$ are in the same component $C$. We may assume $C = C_i$ for some $i > 1$ since $C_1$ is a subgraph of $G[0, b_0]$. Let $e$ be the cut-edge immediately before $C_i$. Suppose that the edge $e$ appears in the segment $G[b_{j-1}, b_j]$. We cannot have $j = 0$, as then we have the contradiction
        \begin{equation*}
            \diam(G_n) = d_{G_n}(u, v) \le a_1 < \frac{c(d)}{2} a_0 \le L[0,b_0] \le d_{G_n}(X_0,u),
        \end{equation*}
        Therefore, $j \ge 1$. Then, we have $d_{G_n}(u, v) \le |E(C_i)| \le a_j + a_{j+1}-1$ since each segment has at least one cut-edge. On the other hand, we have
        \begin{equation*}
            d_{G_n}(X_0, u) \ge \sum_{i'=0}^{j-1} L[b_{i'-1},b_{i'}] \ge \frac{c(d)}{2} \sum_{i'=0}^{j-1} a_{i'} \ge a_j + a_{j+1} > d_{G_n}(u,v) = \diam(G_n),
        \end{equation*}
        which is a contradiction (the third inequality above uses \Cref{lem:cond}). Therefore, $u$ and $v$ belong to different components, say $C_{i_1}$ and $C_{i_2}$, respectively. We may assume $i_1 < i_2$ without loss of generality. Let $x$ be the left endpoint of the cut-edge to the immediate after $C_{i_1}$. We define
        \begin{equation*}
            j = \max \{j' \ge 0: x \in V(G[b_{j'-1}+1,b_{j'}])\}.
        \end{equation*}
        We can assume $j \ge 1$, otherwise $u$ appears in the segment $G[0,b_0]$. If $j = 1$, then we have
        \begin{align*}
            \diam(G_n) = d_{G_n}(u, v) &= d_{G_n}(u, x) + d_{G_n}(x, v)\\
            &\le a_1 + d_{G_n}(x, v)\\
            &< \frac{c(d)}{2} a_0 + d_{G_n}(x, v)\\
            &\le L[0,b_0] + d_{G_n}(x, v)\\
            &\le d_{G_n}(X_0, x) + d_{G_n}(x, v) = d_{G_n}(X_0, v),
        \end{align*}
        a contradiction. Therefore, $j \ge 2$. Then, we have $d_{G_n}(u, x) \le |E(C_{i_1})| \le a_j + a_{j-1} - 1$. On the other hand, we have
        \begin{equation*}
            d_{G_n}(X_0, x) \ge \sum_{i'=0}^{j-1} L[b_{i'-1},b_{i'}]  \ge \frac{c(d)}{2} \sum_{i'=0}^{j-2} a_{i'} \ge a_j + a_{j-1}.
        \end{equation*}
        Using a similar argument as above, we conclude $d_{G_n}(X_0, v) > d_{G_n}(u, v) = \diam(G_n)$, which is a contradiction. Hence, at least one of $u$ and $v$ must appear in the first segment $G[0,b_0]$.
    \end{proof}
    
    Suppose without loss of generality that $u \in V(G[0, b_0])$. Then the source is contained in the closed ball $B_{G_n}(u, m)$ of radius $m$ around $u$ in $G_n$, which is in turn contained inside $B_{G_n}(X_0, 2m)$. Let $i_0 = \min \{i \in \NN: b_i \ge 4m/c(d)\}$. Then the number of cut-edges in the graph $G[0, b_{i_0}]$ is at least $c(d)b_{i_0}/2 \ge 2m$. Thus, it follows that $B_{G_n}(X_0, 2m) \subseteq V(G[0, b_{i_0}])$. Consequently, we have
    \begin{align*}
        |B_{G_n}(X_0, 2m)| &\le |V(G[0, b_{i_0}])| \le b_{i_0} + 1 \le b_{i_0 - 1} + a_{i_0} + 1\\
        &\le b_{i_0 - 1} + a_{i_0} + a_{i_0 + 1} \le b_{i_0 - 1} \left(1 + \frac{c(d)}{2}\right) < 2m \left(\frac{2}{c(d)} + 1\right),
    \end{align*}
    which readily implies the desired result.
\end{proof}

We now show that the event $\cA_m$ occurs with probability approaching $1$ as $m \to \infty$.

\begin{lemma}
\label{lem:pa}
    We have
    \begin{equation*}
        \PP(\cA_m) = 1 - O(m^{-1/2}) \text{ as } m \xrightarrow{} \infty.
    \end{equation*}
\end{lemma}

\begin{proof}
    By \Cref{lem:var}, there is $C(d) > 0$ such that
    \begin{equation*}
        \text{Var}(|L_k|) \le C(d) k^{3/2}, \text{ for all } k \in \NN.
    \end{equation*}
    For each $n \ge 0$, we have
    \begin{equation*}
        \PP\left(|L[b_{n-1}, b_n]| < \frac{c(d)}{2} a_n\right) = \PP\left(|L_{a_n}| < \frac{c(d)}{2} a_n\right) 
        \le \PP\left(||L_{a_n}| - \EE[|L_{a_n}|]| \ge \frac{c(d)}{2} a_n\right).
    \end{equation*}
    Applying Chebychev's inequality to  $|L_{a_n}|$, we conclude
    \begin{equation*}
        \PP\left(|L[b_{n-1}, b_n]| < \frac{c(d)}{2} a_n\right) \le \frac{4 C(d)}{c(d)^2} a_n^{-1/2}.
    \end{equation*}
    Taking a union bound over $n \in \ZZ_{\ge 0}$, we obtain
    \begin{equation*}
        \PP(\cA_m^\complement) \le \frac{4 C(d)}{c(d)^2} \sum_{n=0}^\infty a_n^{-1/2}.
    \end{equation*}
    Now, as the sequence $(a_n)_{n \ge 0}$ grows geometrically by \Cref{lem:cond}, the sum on the right-hand side of the above inequality is $O(m^{-1/2})$.
\end{proof}

We now have all the ingredients to prove the main result of this section.

\begin{proof}[Proof of \Cref{thm:SRWconst}]
    \label{prf:thmSRWconst}
    By \Cref{lem:pa}, there is $C > 0$ such that $\PP(\cA_m) \ge 1 - C/\sqrt{m}$. Let $\varepsilon \in (0, 1)$, choose $m \in \NN$ so that $C/\sqrt{m} < \varepsilon$, and let $k = 2 \lfloor 2m (\frac{2}{c(d)} + 1) \rfloor$. Then \Cref{lem:agood} implies that $\Lambda_k(G_n)$ contains the source whenever $\cA_m$ occurs. Thus, we obtain
    \begin{equation*}
        \PP\left(X_0 \in \Lambda_k(G_n)\right) \ge \PP(\cA_m) \ge 1 - C/\sqrt{m} \ge 1 - \eps. \qedhere
    \end{equation*}
\end{proof}

\subsection{Infinite trace}
\label{sec:infinite}

In this section, we consider the problem of localising the source of the infinite trace of the (one-sided) simple random walk. Our localising method is based on our approach in the finite case.

Let $G_\infty$ be the total trace of the simple random walk on $\ZZ^d$, i.e., $V(G_\infty) = \{X_n: n \in \ZZ_{\ge 0}\}$ and $E(G_\infty) = \{\{X_n, X_{n+1}\}: n \in \ZZ_{\ge 0}\}$. We define an estimator $\Gamma$ as follows.
\begin{itemize}
    \item Let $x$ be an arbitrary vertex in $V(G_\infty)$, chosen in a measurable way.
    \item For $n \ge 1$, consider $B_n = B_G(x,n)$, i.e., the ball of radius $n$ around the vertex~$x$ in~$G_\infty$. Consider the graph obtained by deleting all edges between vertices in $B_n$ from~$G_\infty$. Among its connected components, a.s., there is exactly one which is infinite; call it~$I_n$. Let $H_n = G_\infty \setminus I_n$ be the graph obtained by deleting the edges of $I_n$ from $G_\infty$ and thereafter deleting isolated vertices.
    \item For every diametrical path $D \in \Diam(H_n)$, let $u_D$ and $v_D$ be its endpoints so that $d_{H_n}(x,u_D) \le d_{H_n}(x,v_D)$ (breaking ties arbitrarily). Define ${\bf U}_n^x = \{u_D \colon D \in \Diam(H_n)\}$.
    \item Define the estimator $\Gamma$ on input $G_\infty$ to be a uniformly chosen point from ${\bf U}_N^x$, where $N$ is given by the following proposition.
\end{itemize}

\begin{proposition}
\label{prop:infstab}
    There is an a.s. finite $N = N(G_\infty)$ such that ${\bf U}_n^x = {\bf U}_N^x$ for every $n \ge N$. Further, we have $\inf_{v\in \ZZ^d} \PP_v({\bf U}_N^x = \{v\}) \ge c(d)$.
\end{proposition}

As a straightforward corollary from \Cref{prop:infstab} and \Cref{lem:cd}, we obtain the following main result of this section.

\begin{theorem}
\label{thm:inftrace}
   For each $d \ge 5$, we have $\inf_{v \in \ZZ^d} \PP_v(\Gamma(G_\infty) = v) > 0$. Moreover, we have
   \begin{equation*}
       \inf_{v \in \ZZ^d} \PP_v(\Gamma(G_\infty) = v) = 1 - O(1/d) \ \text{ as } \ d \to \infty.
   \end{equation*}
\end{theorem}

For the proof, we start with an auxiliary lemma (notation is as in the previous section).

\begin{lemma}
\label{lem:sandwich}
    Let $m \ge 8 /c(d)$. There are constants $k_0 = k_0(d) \in \NN$ and $c_0 = c_0(d) > 0$ such that the following holds. Suppose that the event $\cA_m$ occurs and $H$, $K$ are connected graphs such that
    \begin{equation*}
        G[0,b_k] \subseteq H \subseteq G[0,b_{k+2}]\quad \text{and} \quad G[0,b_l] \subseteq K \subseteq G[0,b_{l+2}]
    \end{equation*}
    for some $k,l \ge k_0$. Let $P \in \Diam(H)$. Then 
    \begin{enumerate}
        \item one of the endpoints, $u$, of $P$ is contained in $G[0,b_0]$,
        \item the other endpoint of $P$ is contained in $G[b_{\lfloor kc_0 \rfloor},b_{k+2}]$, and
        \item there is $Q \in \Diam(K)$ which also has $u$ as one of its endpoints.
    \end{enumerate}
\end{lemma}

\begin{proof}
     \textit{(i)}: The proof follows in the same fashion as the proof of \Cref{lem:agood} with the only difference being the use of $H$ in place of $G_n$. We remark that the proof does not use the assumption that $H \subseteq G[0,b_{k+2}]$.

    \textit{(ii)}: Denote the other endpoint of $P$ by $v$. Let $s = \lfloor kc(d) / 3\rfloor - 2$ and suppose that $v \in G[0,b_s]$. Since $u \in G[0,b_0]$, and $G[b_s,b_{s+1}]$ contains some cut-edge, we conclude that $P \subseteq G[0,b_{s+1}]$. Thus, we obtain
    \begin{equation*}
        \diam(H) \le |E(G[0,b_{s+1}])| \le b_{s+1}.
    \end{equation*}

    Further, since $G[0,b_k] \subseteq H$, by \eqref{eq:eventAm}, we get that the total number of cut-edges contained in $H$ is at least $b_kc(d)/2$ and so 
    \begin{equation}\label{eq:diampart1}
        b_{s+1} \ge \diam(H) \ge b_kc(d)/2.
    \end{equation} 

    Observe that for $i \ge 1$, we have $a_{i+1} \ge a_i$. Further, it follows from \Cref{lem:cond} that for $l \ge l(d)$, where $l(d) = 2 + \log(1/\lambda)/\log (1 + \lambda/(1+\lambda))$, we have $a_l \ge a_0$. Thus, if $s+2 \ge l(d)$, we get
    \begin{align*}
        b_k \ge b_{\left\lfloor k / (s+2)\right\rfloor (s+2) - 1} = \sum_{i=0}^{\lfloor \tfrac{k}{s+2}\rfloor-1} \sum_{j=0}^{s+1} a_{i(s+2)+j} \ge \lfloor \tfrac{k}{s+2}\rfloor \sum_{j=0}^{s+1} a_j = \lfloor \tfrac{k}{s+2}\rfloor b_{s+1} \ge \lfloor \tfrac{3}{c(d)}\rfloor b_{s+1},
    \end{align*}
    Finally, $c(d) \le 1$, so $\lfloor \tfrac{3}{c(d)}\rfloor > \tfrac{3}{c(d)} - 1 \ge \tfrac{2}{c(d)}$. Thus we obtain $b_k > 2b_{s+1}/c(d)$, which is a contradiction with \eqref{eq:diampart1}. As a result, $v \in G[b_s,b_{k+2}]$, which implies the statement with $k_0 = \lceil(18/c(d))(l(d)+1)\rceil$ and $c_0 = c(d)/6$.

    \textit{(iii)}: Choose an arbitrary $R \in \Diam(K)$. By the first two parts, its endpoints $x$ and $y$ satisfy $x \in G[0,b_0]$ and $y \in G[b_{\lfloor l c_0 \rfloor}, b_{l+2}]$. As a result, both $R$ and $P$ contain the cut-edges contained in $G[b_0,b_1]$; let $z$ be an endpoint of one of these cut-edges. Clearly, we have $d_G(z,u) = d_G(z,x)$. Substituting the segment of $P$ between $u$ and $z$ for the segment of $R$ between $x$ and $z$ yields a desired diametrical path of $K$.
\end{proof}

\begin{proof}[Proof of \Cref{prop:infstab}]\label{prf:propinfstab}
    Suppose that $\cA_m$ holds for some $m \ge 8/c(d)$. Let $n_0$ be the first time when the walk visits $x$. We begin with the following claim.

    \begin{claim}
    \label{cl:sandwich}
        For every $n > n_0 + b_{k_0}$ there is $k \ge c_1$ such that $G[0,b_k] \subseteq H_n \subseteq G[0,b_{k+2}]$.
    \end{claim}

    \begin{proof}[Proof of \Cref{cl:sandwich}]
        Enumerate the cut-edges by $e_1, e_2, e_3, \ldots$, and let $t_1 < t_2 < t_3 < \ldots$ be the times when they respectively appear (i.e. $e_i = X[t_i-1,t_i]$). We will show that if $e_i \in H_n$, then $G[0,t_i] \subseteq H_n$. Indeed, recall that $H_n = G_\infty \setminus I_n$, where $I_n$ is an infinite connected subgraph of~$G_\infty$. As $e_i$ is a cut-edge, we have $I_n \subseteq G[t_i,+\infty)$, and therefore $G[0,t_i] \subseteq H_n$, as required.
        
        As $n > n_0 + b_1$, we have $H_n \supseteq G[0,b_1]$ and so $H_n$ contains at least one cut-edge (one from $G[b_0,b_1])$. Let then $e_l$ be the last cut-edge contained in $H_n$. We have $t_l \in [b_k,b_{k+1}]$ for some $k \ge 0$. As a result, 
        \begin{equation*}
            G[0,b_k] \subseteq G[0,t_l] \subseteq H_n.
        \end{equation*}
        Further, since $e_l$ is the last cut-edge within $H_n$, we have $t_{l+1} \in [b_k,b_{k+2}]$ and so 
        \begin{equation*}
            H_n \subseteq G[0,t_{l+1}] \subseteq G[0,b_{k+2}].
        \end{equation*}
        The fact that $n > n_0 + b_{k_0}$ enforces $k \ge k_0$, and so the claim is confirmed.
    \end{proof}
    
    Take $N = n_0 + b_{\max(k_0, 2n_0+2m)}$. Consider $i\ge N$. By the above claim, we can apply \Cref{lem:sandwich} to $H = H_i$. Let $D \in \Diam(H_i)$ and denote its endpoints by $y$, $z$. By \Cref{lem:sandwich} (i) and (ii), exactly one of $y$, $z$ is contained in $G[0,b_0]$; without loss of generality, suppose it is $y$.

    We have
    \begin{equation*}
        d_{H_i}(x,y) \le d_{H_i}(x,X_0) + d_{H_i}(X_0,y) \le n_0 + b_0 = n_0+m.
    \end{equation*}
    Observe that as $i \ge N$, there is $H_i \supseteq G[0,b_{2n_0+2m}]$ and therefore $\diam(H_i) > 2(n_0+m)$. Thus
    \begin{equation*}
        d_{H_i}(x,z) \ge d_{H_i}(y,z) - d_{H_i}(x,y) > 2(n_0+m) - (n_0+m) \ge d_{H_i}(x,y).
    \end{equation*}
    Hence we get that $y = u_D$ and $z = v_D$, and so $u_D \in G[0,b_0]$.  

    It follows that for every $i \ge N$ we have ${\bf U}_i^x = {\bf U}_N^x.$ Indeed, take any $D \in \Diam(H_N)$. By the above, $u_D \in G[0,b_0]$ and so (iii) of \Cref{lem:sandwich} applies to it. As a result, there is $D' \in \Diam(H_i)$ such that $u_D$ is its endpoint. As $u_D \in G[0,b_0]$, by the above, we must have $u_{D'} = u_D$. This yields ${\bf U}_N^x \subseteq {\bf U}_i^x$. Analogously, we obtain ${\bf U}_N^x \supseteq {\bf U}_i^x$, and thus ${\bf U}_N^x = {\bf U}_i^x$.

    For the proof of the second part, let $R_m$ be the event that the source is an endpoint of every diametrical path of $G[0, b_{k_0}]$. As in the \hyperref[prf:thmSRW1]{proof} of \Cref{thm:SRW}, recalling that $\cC_n$ is the event that the diameter increases at step $n+1$, by reversibility we have
    \begin{equation}
    \label{eq:rmlower}
        \PP(R_m) \ge \PP(\cC_{b_{k_0}-1}). 
    \end{equation}

    Suppose now that both $R_m$ and $\cA_m$ hold. Let $D$ be any diameter of $H_N$ and consider $u_D$. Applying \Cref{lem:sandwich} to $H_N$ and $G[0,b_{k_0}]$, and diameter $D$ of $H_N$, by (iii) we obtain that there is a diameter of $G[0,b_{k_0}]$ whose endpoint is also $u_D$. By the above, we know that $u_D \in G[0,b_0]$ and so, as the event $R_m$ holds, we obtain $u_D = X_0$. Therefore, whenever the event $\cA_m \cap R_m$ occurs, we have ${\bf U}^x_N = \{X_0\}$. Hence,
    \begin{align*}
        \PP_v({\bf U}^x_N = \{v\}) &\ge  \limsup_{m \to \infty} \PP(\cA_m \cap R_m) \ge \limsup_{m\to\infty} \ (\PP(R_m) - \PP(\cA_m^\complement) ) \ge c(d),
    \end{align*}
    where the last inequality follows by \Cref{lem:pa}, \eqref{eq:rmlower} and \Cref{prop:PC}.
\end{proof}

The above yields the proof of \Cref{thm:inftrace}, which is a restatement of the first part of \Cref{thm:infinite}. Now we turn to proving the second part of \Cref{thm:infinite}.

\begin{proof}[Proof of the upper bound in \Cref{thm:infinite}]
\label{prf:upperboundinfinite}
    The proof is similar to and simpler than the \hyperref[prf:ub]{proof} of the upper bound in \Cref{thm:SRWoptimal}. Let $E \coloneqq \{\bfe_1, \dots, \bfe_d\}$ be the set of standard basis vectors for $\ZZ^d$, and let $E_\pm \coloneqq \{\pm \bfe: \bfe \in E\}$. Let $\Omega \coloneqq E_\pm^\NN$ be the space of sequences with entries in $E_\pm$. We define the left shift map $f: \Omega \to \Omega$ by
    \begin{equation*}
        f((s_1, s_2, s_3, s_4, \dots)) = (s_2, s_3, s_4, \dots).
    \end{equation*}
    Let $S = (S_1, S_2, S_3, S_4, \dots)$ be a sequence of iid random variables with uniform distribution on $E_\pm$. Note that $f(S)$ has the same law as $S$, namely, the product measure $\mu = \otimes_{j \ge 1} \nu$, where $\nu$ is the uniform measure on $E_\pm$. Let $X$ and $X^f$ be simple random walks on $\ZZ^d$ started at the origin with steps given by $S$ and $f(S)$, respectively. Further, let us denote the corresponding infinite trace graphs by $G_\infty$ and $G^f_\infty$, respectively.
    
    For $\ell \in \NN$, let us define the box $B_\ell \coloneqq [-\ell, \ell]^d \cap \ZZ^d$ and consider the random variable
    \begin{equation*}
        W \coloneqq \sum_{\bfv \in B_\ell} \mathbf{1}\{\Phi(G_\infty + \bfv) = \bfv\} + \sum_{\bfv \in B_\ell} \mathbf{1}\{\Phi(G^f_\infty + \bfv + X_1) = \bfv + X_1\}.
    \end{equation*}
    Taking expectation on both sides, we obtain
    \begin{equation}
    \label{eqn:expW}
        \EE[W] = \sum_{\bfv \in B_\ell} \PP(\Phi(G_\infty + \bfv) = \bfv) + \sum_{\bfv \in B_\ell} \PP(\Phi(G^f_\infty + \bfv + X_1) = \bfv + X_1).
    \end{equation}
    By a similar argument as in the \hyperref[prf:ub]{proof} of the upper bound in \Cref{thm:SRWoptimal}, the second sum on the right-hand side above is bounded below by
    \begin{equation*}
        \sum_{\bfv' \in B_{\ell-1}} \PP(\Phi(G_\infty^f + \bfv') = \bfv').
    \end{equation*}
    Using the above estimate in \eqref{eqn:expW}, we obtain
    \begin{align}
    \label{eqn:EWlow}
        \EE[W] &\ge \sum_{\bfv \in B_\ell} \PP(\Phi(G_\infty + \bfv) = \bfv) + \sum_{\bfv' \in B_{\ell-1}} \PP(\Phi(G_\infty^f + \bfv') = \bfv')\nonumber \\
        &\ge (|B_\ell| + |B_{\ell-1}|) \inf_{\bfv \in \ZZ^d} \PP(\Phi(G_\infty + \bfv) = \bfv).
    \end{align}
    
    On the other hand, let $A$ be the event $\{-S_1 = S_2\}$. Note that $X_1 \ne 0$ and whenever $A$ occurs, we have $G_\infty^f + X_1 = G_\infty$. Thus, whenever $A$ occurs, for each $v \in \ZZ^d$, we have
    \begin{equation*}
        \mathbf{1}\{\Phi(G_\infty + \bfv) = \bfv\} + \mathbf{1}\{\Phi(G^f_\infty + \bfv + X_1) = \bfv + X_1\} \le 1.
    \end{equation*} 
    Summing the above over $\bfv \in B_\ell$ implies that $W \le |B_\ell|$ whenever $A$ occurs. Moreover, it follows from the definition of $W$ that $W \le 2 |B_\ell|$. Hence, we obtain
    \begin{equation*}
        \EE[W] = \PP(A) \EE[W|A] + \PP(A^\complement) \EE[W|A^\complement] \le \PP(A)|B_\ell|+ (1 - \PP(A)) \cdot 2 |B_\ell| = (2 - \PP(A))|B_\ell|.
    \end{equation*}
    Combining the above estimate with \eqref{eqn:EWlow}, we obtain
    \begin{equation*}
        (|B_\ell| + |B_{\ell-1}|) \inf_{\bfv \in \ZZ^d} \PP(\Phi(G_\infty + \bfv) = \bfv) \le (2 - \PP(A))|B_\ell|.
    \end{equation*}
    Observe that $|B_\ell| = \Theta(\ell^d)$ and $|B_{\ell-1}| = |B_\ell| - O(\ell^{d-1})$ as $\ell \rightarrow \infty$. Therefore, dividing both sides of the above inequality by $|B_\ell| + |B_{\ell-1}|$, and taking $\lim_{\ell \to \infty}$, we obtain
    \begin{equation*}
        \inf_{\bfv \in \ZZ^d} \PP(\Phi(G_\infty + \bfv) = \bfv) \le 1 - \frac{\PP(A)}{2}.
    \end{equation*}
    Substituting $\PP(A) = 1/(2d)$ yields the desired result.
\end{proof}

\subsection{Range input}
\label{sec:range}

In this section, we consider the variant where the input is the trace consisting of exactly $n$ vertices as opposed to corresponding to the first $n$ steps of the random walk. Let
\begin{equation*}
    \tau_i = \inf \{t\ge 0\colon |V(G_t^d)| = i\},
\end{equation*}
i.e., $\tau_i$ is the first time the trace graph consists of $i$ vertices. Define $R_i^d \coloneqq G^d_{\tau_i}$. Recall that the estimator $\Psi$ samples uniformly an endpoint of a uniformly selected diametrical path of $G$.

\begin{proof}[Proof of \Cref{thm:range}]
    Repeating the argument of the second part of the \hyperref[prf:propinfstab]{proof} of \Cref{prop:infstab}, we obtain that there is an event $\cE_n$ such that $\PP(\cE_n) \to c(d)$ as $n\to \infty$ and whenever $\cE_n$ occurs, for $m \ge n$, every diametrical path of $G[0,m]$ has $X_0$ as an endpoint. 
    
    Further, observe that $\tau_n \ge n$ a.s. Hence, whenever $\cE_n$ occurs, every diameter of $R_n^d = G_{\tau_n}^d$ has $X_0$ as an endpoint. Thus,
    \begin{equation*}
        \PP( \Psi(R_n^d)=X_0) \ge \PP(\cE_n) \cdot \PP(\Psi(G_{\tau_n}^d) = X_0 \mid \cE_n) = \PP(\cE_n)/2.
    \end{equation*}
    Taking $\liminf_{n\to\infty}$ of both sides yields the statement.
\end{proof}

\subsection{Vertex trace only}
\label{sec:vertex}

Let $(\overline{X}_n)_{n \in \ZZ}$ be the two-sided random walk on $\ZZ^d$. Recall that for $n \in \ZZ$, we defined $I_n$ as the indicator function of the event $\left\{\overline{X}(-\infty, n] \cap \overline{X}[n+1,\infty) = \emptyset\right\}$. In a similar way, we define $J_n$ as the indicator of the event
\begin{equation*}
    \left\{\left\|\overline{X}_i - \overline{X}_j\right\|_1 \ne 1 \text{ for all }i \le n, j \ge n+1 \text{ with } (i,j) \ne (n, n+1)\right\}.
\end{equation*}
If $J_n = 1$, we say that the edge $\{X_n, X_{n+1}\}$ is an \emph{induced cut-edge} of the two-sided random walk. This notion has previously been studied under the name of strong cut points \cite{JLP}. Let $V$ denote the vertex set $\overline{X}(-\infty, \infty)$. Similarly to cut-edges, the key observation is that induced cut-edges of the two-sided random walk are cut-edges of the induced graph $\ZZ^d[V]$. Then, the proof of \Cref{thm:SRW2} works the same way as that of \Cref{thm:SRW} except that we now use induced cut-edges instead of cut-edges. Similarly to the definition of $c(d)$, we define
\begin{equation*}
    \widetilde{c}(d) \coloneqq \EE[J_0] = \PP(J_0 = 1).
\end{equation*}
To finish the proof of \Cref{thm:SRW2}, it suffices that $\widetilde{c}(d) > 0$ for $d \ge 5$ and $\widetilde{c}(d) = 1 - O\left(1/d\right)$. 

These estimates are fairly standard, but for completeness, we include the proofs below. Let $(X_n)_{n \ge 0}$ and $(X'_n)_{n \ge 0}$ be simple random walks started at adjacent vertices in $\ZZ^d$, say $X_0 = 0$ and $X'_0 = (1, 0, \dots, 0)$. Let $J$ be the number of times the two random walk traces are adjacent in $\ZZ^d$. More precisely, we define
\begin{equation*}
    J \coloneqq \left|\left\{(i, j) \in \ZZ_{\ge 0}^2 \colon \left\|X_i - X'_j\right\|_1 = 1\right\}\right|.
\end{equation*}
In other words, $J$ counts the number of \emph{adjacencies} between the two random walks.

\begin{lemma}
\label{lem:EI}
    For each $d \ge 5$, $\EE[J]$ is finite. Moreover, we have
    \begin{equation*}
        \EE[J] = 1 + O(1/d).
    \end{equation*}
\end{lemma}

\begin{proof}
    We have
    \begin{align*}
        \EE[J] &= \sum_{i = 0}^\infty \sum_{j = 0}^\infty \PP\left(\left\|X_i - X'_j\right\|_1 = 1\right)
        = \sum_{i = 0}^\infty \sum_{j = 0}^\infty \PP\left(\left\|X_{i+j+1}\right\|_1 = 1\right)\\
        &= \sum_{n = 0}^\infty (2n+1) \PP\left(\left\|X_{2n+1}\right\|_1 = 1\right) = 1 + 2d \sum_{n = 1}^\infty (2n+1) \PP(X_{2n+2} = 0)\\
        &\le 1 + 4d \sum_{n = 2}^\infty n \PP(X_{2n} = 0).
    \end{align*}
    The result follows from \Cref{lem:npsum} with $k = 2$.
\end{proof}

We relate $\widetilde{c}(d)$ and $\EE[J]$ in the following lemma, which is similar to \cite[Proposition~3.2.2]{Law2}.

\begin{lemma}
\label{lem:cdE}
    For each $d \ge 5$, we have
    \begin{equation*}
        \widetilde{c}(d) \ge \EE[J]^{-1}.
    \end{equation*}
    Consequently, $\widetilde{c}(d) > 0$ for each $d \ge 5$, and
    \begin{equation*}
        \widetilde{c}(d) = 1 - O(1/d).
    \end{equation*}
\end{lemma}

\begin{proof}
    We say that $(i, j) \in \ZZ_{\ge 0}^2$ is a \emph{maximal adjacency} if $\left\|X_i - X'_j\right\|_1 = 1$ and
    \begin{equation*}
        \left\|X_{i'} - X'_{j'}\right\|_1 \ne 1, \text{ for all } i' \ge i,\, j' \ge j,\, (i', j') \ne (i, j).
    \end{equation*}
    Let $M$ be the number of maximal adjacencies. Then, we have
    \begin{align*}
        \EE[M] &= \sum_{i = 0}^\infty \sum_{j = 0}^\infty \PP((i,j) \text{ is a maximal adjacency})\\
        &= \sum_{i = 0}^\infty \sum_{j = 0}^\infty \PP\left(\left\|X_i - X'_j\right\|_1 = 1\right) \cdot \PP((X_n)_{n \ge i} \text{ and } (X'_n)_{n \ge j} \text{ have no non-trivial adjacency})\\
        &= \widetilde{c}(d) \EE[J].
    \end{align*}
    Note that $M \ge 1$ almost surely since the set of adjacencies is non-empty and almost surely finite (see the first part of \Cref{lem:EI}). This implies $\EE[M] \ge 1$, which when used in the above equation yields the first bound. The second bound follows then from \Cref{lem:EI}.
\end{proof}

\section{Concluding remarks} \label{sec:conclude}

By \Cref{thm:SRW} and \Cref{thm:SRWrecurrent}, the simple random walk on a locally finite vertex-transitive graph is not amnesiac when the graph is strongly transient, while it is amnesiac when the graph is recurrent. In particular, the simple random walk on $\ZZ^d$ is not amnesiac for $d \ge 5$, while it is amnesiac for $d \le 2$. This leaves the following question for the borderline dimensions $d = 3, 4$.

\begin{question}
\label{ques:34}
    Are the simple random walks on $\ZZ^3$ and $\ZZ^4$ amnesiac?
\end{question}

More generally, for random walks on discrete groups, we note that the group must have a polynomial growth of degree $3$ or $4$ for the simple random walk on it to be transient but not strongly transient; see \cite[page~622]{Bla}. In this case, it turns out that the group must contain $\ZZ^3$, $\ZZ^4$, or the $3$-dimensional discrete Heisenberg group $H_3(\ZZ)$ as a finite index subgroup \cite[Proposition~2.5]{Bla}. Therefore, it would also be interesting to study \Cref{ques:34} for $H_3(\ZZ)$. In fact, we are not aware of any examples of vertex-transitive transient graphs on which the simple random walk is amnesiac. This suggests the following broad question.

\begin{question}
    Is there any locally finite vertex-transitive transient graph on which the simple random walk is amnesiac? In other words, is amnesia equivalent to recurrence for such graphs?
\end{question}

We remark that any locally finite vertex-transitive transient graph on which the simple random walk is amnesiac (if there exists one) must exhibit polynomial growth with growth exponent less than $5$ since otherwise the walk is known to be strongly transient~\cite[Corollary~14.5]{Woe}.

Regarding the high accuracy localisation problem, the \hyperref[prf:thmSRWconst]{proof} of \Cref{thm:SRWconst} shows that one can take $K(\varepsilon) = O(\varepsilon^{-2})$, with a better exponent for large $d$. We do not make the dependence explicit, as we have not investigated the sharp dependence, which may be considered in future research. We do at least note that for some $n=n(\eps)$ one needs $K(\varepsilon)$ to grow polynomially in $1/\varepsilon$, as may be seen by considering the event that the final vertex $X_n$ is equal to the source vertex $X_0$, in which case the walk becomes invariant under cyclic shifts of the source. On the other hand, we would not be surprised if logarithmic growth in $1/\varepsilon$ suffices for very large $n$.

\bibliographystyle{amsplain}
\bibliography{bibliography}

\appendix

\section{Estimates for the probability of a cut-edge}
\label{app:cd}

In the following, our main objective is to obtain estimates on the growth of the cut-edge probability $c(d)$ for $\ZZ^d$ as $d \rightarrow \infty$. We will require a rough estimate for the probability of the simple random walk returning to the origin after a given number of steps.

\begin{lemma}[Return probability estimate]
\label{lem:lclt}
    Let $d \ge 1$ and let $(X_n)_{n \ge 0}$ be a simple random walk on $\ZZ^d$. Then, we have
    \begin{equation*}
        \PP(X_{2n} = 0) \le \left(\frac{2d}{\pi n}\right)^{\frac{d}{2}}, \text{ for all } n \in \NN.
    \end{equation*}
\end{lemma}

\begin{proof}
    Let $\phi(\theta) \coloneqq \frac{1}{d} \sum_{j=1}^d \cos \theta_j$ for all $\theta \in \RR^d$,
    which is the characteristic function of $X_1$. Then  $\phi^n$ is the characteristic function of $X_n$. 
    By the Fourier inversion formula, we have
    \begin{align}
    \label{eqn:est}
        \PP(X_{2n} = 0) & = (2\pi)^{-d} \int_{[-\pi,\pi]^d} \phi^{2n}(\theta) \d\theta, \nonumber \\
        &\le (2\pi)^{-d} \int_{[-\pi,\pi]^d} \left(\frac{1}{d} \sum_{j=1}^d \left|\cos \theta_j\right|\right)^{2n} \d\theta \nonumber\\
        &= \frac{1}{\pi^d} \int_{[-\pi/2, \pi/2]^d} \phi^{2n}(\theta) \d\theta.
    \end{align}
   Using $0 \le \cos \varphi \le 1 - \varphi^2/4$ for $|\varphi| \le \pi/2$ we have
    \begin{equation*}
        0 \le \phi(\theta) \le 1 - \frac{1}{4d}\|\theta\|_2^2 \le \e^{-\frac{\|\theta\|_2^2}{4d}}, \text{ for all } \theta \in [-\pi/2, \pi/2]^d.
    \end{equation*}
    Using this in \eqref{eqn:est} and recognising the density integral for a $d$-variate normal yields
    \begin{equation*}
        \PP(X_{2n} = 0) \le \frac{1}{\pi^d} \int_{[-\pi/2, \pi/2]^d} \e^{-\frac{n\|\theta\|_2^2}{2d}} \d\theta \le \frac{1}{\pi^d} \int_{\RR^d} \e^{-\frac{n\|\theta\|_2^2}{2d}} \d\theta = \left(\frac{2d}{\pi n}\right)^{d/2}. \qedhere
    \end{equation*}
\end{proof}

We will also require the fact that the return probability to the origin for the simple random walk is non-increasing in the number of steps. Pittet and Saloff-Coste~\cite{PS} remark that this is well known in the more general setting of random walks on Cayley graphs with a symmetric generating set. We state this in \Cref{lem:monotonicity} for the simple random walk on $\ZZ^d$. We include a proof for completeness, which also works in more general setting.

\begin{lemma}[Monotonicity of return probability]
\label{lem:monotonicity}
    Let $d \ge 1$ and let $(X_n)_{n \ge 0}$ be a simple random walk on $\ZZ^d$. Then $\PP(X_{2n} = 0)$ is non-increasing in $n$ for $n \in \ZZ_{\ge 0}$.
\end{lemma}

\begin{proof}
    Observe that for any $x \in \ZZ^d$ and $n \ge 0$, we have
    \begin{align*}
        \PP(X_{2n} = x) &= \sum_{y \in \ZZ^d} \PP(X_n = y) \PP(X_n = x - y) 
        \le \sum_{y \in \ZZ^d} \frac{\PP(X_n = y)^2+  \PP(X_n = x - y)^2}{2}\\
        &= \sum_{y \in \ZZ^d} \PP(X_n = y)^2 
        = \sum_{y \in \ZZ^d} \PP(X_n = y)\PP(X_n = -y) 
        = \PP(X_{2n} = 0).
    \end{align*}
    As a result, for each $n \ge 0$, we have
    \begin{align*}
        \PP(X_{2(n+1)} = 0) &= \sum_{y \in \ZZ^d} \PP(X_{2n} = y) \PP(X_2 = -y)\\
        &\le \PP(X_{2n} = 0) \sum_{y \in \ZZ^d} \PP(X_2 = -y)\\
        &= \PP(X_{2n} = 0). \qedhere
    \end{align*}
\end{proof}

Let $(X_n)_{n \ge 0}$ be a simple random walk in $\ZZ^d$. For each $k \in \NN$, we define
\begin{equation}
\label{eqn:sumdef}
    S_d(k) \coloneqq \sum_{n = k}^\infty n \PP(X_{2n} = 0).
\end{equation}

\begin{lemma}
\label{lem:npsum}
    Let $d \ge 5$. Then, $S_d(k)$ is finite for each $k \in \NN$. Moreover, for any fixed $k \in \NN$, we have $S_d(k) = O(d^{-k})$ as $d \rightarrow \infty$.
\end{lemma}

\begin{proof}
    The finiteness follows by using \Cref{lem:lclt} in \eqref{eqn:sumdef} and observing that the sequence $(n^{-d/2+1})_{n \in \NN}$ is summable for $d \ge 5$. For the second part, let $K \ge 2$ be an arbitrary integer constant. We break the sum on the right-hand side of \eqref{eqn:sumdef} into two parts:
    \begin{equation}
    \label{eqn:parts}
        \sum_{n = k}^\infty n \PP(X_{2n} = 0) = \sum_{n = k}^{d^K} n \PP(X_{2n} = 0) + \sum_{n = d^K+1}^\infty n \PP(X_{2n} = 0).
    \end{equation}
    
    We will first bound the first sum on the right-hand side of \eqref{eqn:parts} using explicit calculations and the monotonicity of the return probability. A direct calculation yields
    \begin{align}
    \label{eqn:ret}
        \PP(X_{2n} = 0) &= \frac{1}{(2d)^{2n}} \sum_{\substack{x_1, \ldots, x_d \ge 0 \\ x_1 + \ldots + x_d = n}} \binom{2n}{x_1, x_1, \ldots, x_d,x_d}\nonumber \\
        &= \frac{\binom{2n}{n}}{(2d)^{2n}} \sum_{\substack{x_1, \ldots, x_d \ge 0 \\ x_1 + \ldots + x_d = n}} \binom{n}{x_1, \ldots, x_d}^2\nonumber \\
        &\le \frac{\binom{2n}{n} \binom{n+d-1}{d-1} (n!)^2}{(2d)^{2n}}.
    \end{align}
    Now, note that
    \begin{equation}
    \label{eqn:b1}
        \frac{n!}{d^n} \binom{n+d-1}{d-1} = \prod_{i=1}^{n-1} \left(1 + \frac{i}{d}\right) \le \exp\left(\frac{n^2}{2d}\right).
    \end{equation}
    Further, we have
    \begin{equation}
    \label{eqn:b2}
        \frac{1}{2^{2n}} \binom{2n}{n} \le \frac{1}{\sqrt{n}},
    \end{equation}
    by a standard bound on the central binomial coefficient, and
    \begin{equation}
    \label{eqn:b3}
        n! \le n^n \sqrt{n} \exp(-n+1),
    \end{equation}
    by Stirling's bound. Using \eqref{eqn:b1}, \eqref{eqn:b2}, and \eqref{eqn:b3} in \eqref{eqn:ret} yields
    \begin{equation*}
        \PP(X_{2n} = 0) \le n^n d^{-n} \exp\left(\frac{n^2}{2d} - n + 1\right) \le n^n d^{-n} \exp\left(-\frac{n}{2} + 1\right), \text{ for } n \le d.
    \end{equation*}
    Now suppose $d \ge 2K+k$. Using the above inequality, we obtain
    \begin{equation}
    \label{eqn:c1}
        \sum_{n = k}^{2K+k-1} n \PP(X_{2n} = 0) \le \sum_{n=k}^{2K+k-1} n^{n+1} d^{-n} \exp\left(-\frac{n}{2}+1\right) = O\left(\frac{1}{d^k}\right).
    \end{equation}
    Further, it follows from \Cref{lem:monotonicity} that
    \begin{equation}
    \label{eqn:c2}
        \sum_{n = 2K+k}^{d^K} n \PP(X_{2n} = 0) \le d^{2K} \PP(X_{2(2K+k)} = 0) = O\left(\frac{1}{d^k}\right).
    \end{equation}
    Adding \eqref{eqn:c1} and \eqref{eqn:c2}, we obtain
    \begin{equation}
    \label{eqn:p1}
        \sum_{n=k}^{d^K} n \PP(X_{2n} = 0) = O\left(\frac{1}{d^k}\right).
    \end{equation}

   For the second sum on the right-hand side of \eqref{eqn:parts}, by \Cref{lem:lclt} we have
    \begin{align}
    \label{eqn:p2}
        \sum_{n = d^K+1}^\infty n \PP(X_{2n} = 0) \le \sum_{n = d^K+1}^\infty \left(\frac{2d}{\pi}\right)^{\frac{d}{2}} n^{1-\frac{d}{2}} \le d^{\frac{d}{2}} (d^K)^{2-\frac{d}{2}} = d^{2K - (K-1)d/2} = O\left(\frac{1}{d^k}\right).
    \end{align}
    
    Finally, adding \eqref{eqn:p1} and \eqref{eqn:p2} yields the desired result.
\end{proof}

To estimate $c(d)$, we compare with the following closely related quantity.
Consider two simple random walks $(X_n)_{n \ge 0}$ 
and $(X'_n)_{n \ge 0}$ in $\ZZ^d$ started at the origin. 
Let $s(d)$ be the probability that they do not intersect except at their initial points, i.e.
\begin{equation*}
    s(d) = \PP( \forall (i,j) \in \ZZ^2_{\ge 0} \setminus \{(0,0)\}\colon X_i \neq X'_j).
\end{equation*}
Clearly $s(d)$ is non-decreasing in $d$, by projection on $\ZZ^{d-1}$. Similarly $c(d)$ is non-decreasing in $d$. Erd\H{o}s and Taylor~\cite[Remark~3]{ET} remarked that $s(d) \rightarrow 1$ as $d \rightarrow \infty$. They did not supply a proof, so for completeness, we provide it below.

\begin{lemma}
\label{lem:pd}
    For each $d \ge 5$, we have $s(d) > 0$. Further, we have
    \begin{equation*}
        s(d) = 1 - O(1/d).
    \end{equation*}
\end{lemma}

\begin{proof}
    Let $d \ge 5$. Let $(X_n)_{n \ge 0}$ and $(Y_n)_{n \ge 0}$ be two independent simple random walks in $\ZZ^d$ started at the origin. Further, let $I$ be the number of intersections of the two random walks, including the starting point, namely,
    \begin{equation*}
        I = \left|\left\{(i,j) \in \ZZ_{\ge 0}^2: X_i = Y_j\right\}\right|.
    \end{equation*}
    Then, \cite[Proposition~3.2.2]{Law2} implies
    \begin{equation}
    \label{eqn:rel}
        s(d) \ge \EE[I]^{-1}.
    \end{equation}
    Further, \cite[Proposition~3.2.1]{Law2} implies that $\EE[I]$ is a finite positive constant that depends on $d$, which when used in \eqref{eqn:rel} implies $s(d) > 0$.
    
    For the second part, it suffices to show that $\EE[I] = 1 + O(1/d)$. We have
    \begin{align*}
        \EE[I] = \sum_{i = 0}^\infty \sum_{j=0}^\infty \PP(X_i = Y_j)
        = \sum_{n = 0}^\infty (2n+1) \PP(X_{2n} = 0)
        \le 1 + 3 \sum_{n=1}^\infty n \PP(X_{2n} = 0).
    \end{align*}
    Finally, using \Cref{lem:npsum} with $k = 1$ in the above inequality yields the desired result.
\end{proof}

We now transfer the bound on $s(d)$ given by \Cref{lem:pd} to a bound on $c(d)$ via comparison.

\begin{lemma}
\label{lem:cd}
    For each $d \ge 5$, we have $c(d) > 0$. Further, we have
    \begin{equation*}
        c(d) = 1 - O(1/d).
    \end{equation*}
\end{lemma}

\begin{proof}
    By \Cref{lem:pd}, it suffices to show that $c(d) \ge s(d)$. For each $d \in \NN$, we have
    \begin{align*}
        1 - c(d) &= \PP(\overline{X}(-\infty,0] \cap \overline{X}[1,\infty) \ne \emptyset)\\
        &= \PP(\overline{X}_i = \overline{X}_j \text{ for some integers } i \le 0, j \ge 1)\\
        &\le \PP(\overline{X}_i = \overline{X}_j \text{ for some integers } i \le 0, j \ge 0 \text{ such that } (i,j) \ne (0, 0))\\
        &= 1 - s(d). \qedhere
    \end{align*}
\end{proof}

\section{Strongly transient walks}
\label{app:strongtrans}

Let $\gH$ be an infinite vertex-transitive graph of finite degree. Suppose the simple random walk $(X_n)_{n \ge 0}$ on $\gH$ is transient. Then it is said to be \emph{strongly transient} if
\begin{equation}
\label{eqn:st}
    \EE[\tau_v^+ \mid X_0 = v, \tau_v^+ < \infty] < \infty,
\end{equation}
where $\tau_v^+ \coloneqq \inf \{t > 0: X_t = v\}$ is the first return time to the starting vertex $v$. Let
\begin{equation*}
    p_n = \PP(\tau_v^+ = n \mid X_0 = v) \text{ and } q_n = \PP(X_n = v \mid X_0 = v),
\end{equation*}
for all $n \in \ZZ_{\ge 0}$. The expectation in \eqref{eqn:st} can be rewritten as
\begin{equation*}
    \EE[\tau_v^+ \mid X_0 = v, \tau_v^+ < \infty] = \sum_{n = 0}^\infty n \PP(\tau_v^+ = n \mid X_0 = v, \tau_v^+ < \infty) = \frac{\sum_{n = 0}^\infty n p_n}{\PP(\tau_v^+ < \infty \mid X_0 = v)}.
\end{equation*}
Note that the denominator of the last expression is positive. Thus the walk is strongly transient if and only if $\sum_{n=0}^\infty n p_n < \infty$.

\begin{lemma}
\label{lem:st}
    The random walk $(X_n)_{n \ge 0}$ is strongly transient if and only if $\sum_{n = 0}^\infty n q_n < \infty$.
\end{lemma}

\begin{proof}
    Suppose that $\sum_{n = 0}^\infty n q_n < \infty$. Then we have $\sum_{n = 0}^\infty n p_n < \infty$ since $0 \le p_n \le q_n$ for all $n \in \ZZ_{\ge 0}$. This implies that the walk is strongly transient.

    Conversely, suppose the walk is strongly transient. Consider the two power series
    \begin{equation*}
        f(x) = \sum_{n = 0}^\infty p_n x^n \quad\text{ and }\quad g(x) = \sum_{n = 0}^\infty q_n x^n.
    \end{equation*}
    Their radii of convergence are both at least $1$, as their coefficients are probabilities, so in $[0,1]$. By transience, we also know that the first power series is convergent at $x = 1$. Furthermore, the definitions of $(p_n)_{n \ge 0}$ and $(q_n)_{n \ge 0}$ yield the relation
    \begin{equation*}
        q_n = \sum_{k = 0}^n p_k q_{n - k}
    \end{equation*}
    for all positive integers $n$. By the above and the fact that $q_0 = 1$, for $|x| < 1$, we have
    \begin{equation*}
        g(x) = \frac{1}{1 - f(x)}.
    \end{equation*}
    Now consider the power series
    \begin{equation*}
        f_1(x) = \sum_{n = 0}^\infty n p_n x^n \quad \text{ and } \quad g_1(x) = \sum_{n = 0}^\infty n q_n x^n.
    \end{equation*}
    Then $f_1(x) = f^\prime(x)$ and $g_1(x) = g^\prime(x) = f'(x)/(1-f(x))^2$ for $|x| < 1$. As the walk is strongly transient, the power series for $f_1(x)$ is convergent at $x = 1$. Then by Abel's theorem, we have
    \begin{equation}
    \label{eqn:f1conv}
        \lim_{x \rightarrow 1^-} f^\prime(x) = f_1(1).
    \end{equation}
    For any positive integer $N$, we have
    \begin{equation*}
         \sum_{n = 0}^N n q_n = \lim_{x \rightarrow 1^-} \sum_{n = 0}^N n q_n x^n \le \lim_{x \rightarrow 1^-} \sum_{n = 0}^\infty n q_n x^n 
         = \lim_{x \to 1^-} \frac{f^\prime(x)}{(1 - f(x))^2} = \frac{f_1(1)}{(1 - f(1))^2},
    \end{equation*}
    where the last equality follows from \eqref{eqn:f1conv} and Abel's theorem applied to the power series defining $f(x)$. Thus the power series for $g_1(x)$ converges at $x = 1$, giving $\sum_{n = 0}^\infty n q_n < \infty$.
\end{proof}

The following is immediate from \Cref{lem:st} and \Cref{lem:npsum}.

\begin{proposition}
\label{prop:zdst}
    For each $d \ge 5$, the graph $\ZZ^d$ is strongly transient.
\end{proposition}

Let $(X_n)_{n \ge 0}$ and $(X^\prime_n)_{n \ge 0}$ be simple random walks started at the same vertex in $\gH$. Let $J$ be the number of times the two random walk traces intersect in $\gH$. More precisely, we define
\begin{equation*}
    J \coloneqq \left|\left\{(i, j) \in \ZZ_{\ge 0}^2 \colon X_i = X^\prime_j\right\}\right|.
\end{equation*}
In other words, $J$ counts the number of intersections between the two random walks. 

\begin{lemma}
\label{lem:cEJ}
    Let $\gH$ be an infinite vertex-transitive strongly transient graph. Then $\EE[J] < \infty$.
\end{lemma}

\begin{proof}
    We have
    \begin{align*}
        \EE[J] = \sum_{i = 0}^\infty \sum_{j = 0}^\infty \PP\left(X_i = X^\prime_j\right)
        = \sum_{i = 0}^\infty \sum_{j = 0}^\infty \PP\left(X_{i+j+1} = X_0\right)
        = \sum_{n = 0}^\infty n \PP\left(X_n = X_0\right),
    \end{align*}
    which is finite by \Cref{lem:st}.
\end{proof}

\begin{lemma}
\label{lem:cH}
    Let $\gH$ be an infinite vertex-transitive strongly transient graph. Then
    \begin{equation*}
        \cgH \ge \EE[J]^{-1}.
    \end{equation*}
\end{lemma}

We omit the proof of \Cref{lem:cH} since it is essentially identical to the proof of \cite[Proposition~3.2.2]{Law2} and also similar to the proof of \Cref{lem:cdE}. We conclude with the following immediate consequence of \Cref{lem:cEJ} and \Cref{lem:cH}.

\begin{lemma}
\label{lem:intprob}
    Let $\gH$ be an infinite vertex-transitive strongly transient graph of finite degree. Then the constant $c(\gH)$ as defined in \eqref{eqn:cH} is positive.
\end{lemma}

\end{document}